\documentclass[10pt]{amsart}
\usepackage{latexsym,amsmath,amssymb,amscd}
\def\today{\ifcase \month \or
   January \or February \or March \or April \or
   May \or June \or July \or August \or
   September \or October \or November \or December \fi
   \space\number\day , \number\year}

\begin{document}

\DeclareRobustCommand{\SkipTocEntry}[4]{} 

\makeatletter
\@addtoreset{figure}{section}
\def\thefigure{\thesection.\@arabic\c@figure}
\def\fps@figure{h,t}
\@addtoreset{table}{bsection}

\def\thetable{\thesection.\@arabic\c@table}
\def\fps@table{h, t}
\@addtoreset{equation}{section}
\def\theequation{
\arabic{equation}}
\makeatother

\newcommand{\bfi}{\bfseries\itshape}

\newtheorem{theorem}{Theorem}
\newtheorem{acknowledgment}[theorem]{Acknowledgment}
\newtheorem{algorithm}[theorem]{Algorithm}
\newtheorem{axiom}[theorem]{Axiom}
\newtheorem{case}[theorem]{Case}
\newtheorem{claim}[theorem]{Claim}
\newtheorem{conclusion}[theorem]{Conclusion}
\newtheorem{condition}[theorem]{Condition}
\newtheorem{conjecture}[theorem]{Conjecture}
\newtheorem{construction}[theorem]{Construction}
\newtheorem{corollary}[theorem]{Corollary}
\newtheorem{criterion}[theorem]{Criterion}
\newtheorem{definition}[theorem]{Definition}
\newtheorem{example}[theorem]{Example}
\newtheorem{lemma}[theorem]{Lemma}
\newtheorem{notation}[theorem]{Notation}
\newtheorem{problem}[theorem]{Problem}
\newtheorem{proposition}[theorem]{Proposition}
\newtheorem{remark}[theorem]{Remark}
\numberwithin{theorem}{section}
\numberwithin{equation}{section}

\newcommand{\todo}[1]{\vspace{5 mm}\par \noindent
\framebox{\begin{minipage}[c]{0.95 \textwidth}
\tt #1 \end{minipage}}\vspace{5 mm}\par}

\newcommand{\1}{{\bf 1}}
\newcommand{\Ad}{{\rm Ad}}
\newcommand{\Alg}{{\rm Alg}\,}
\newcommand{\Aut}{{\rm Aut}\,}
\newcommand{\ad}{{\rm ad}}
\newcommand{\Ci}{{\mathcal C}^\infty}
\newcommand{\de}{{\rm d}}
\newcommand{\ee}{{\rm e}}
\newcommand{\ev}{{\rm ev}}
\newcommand{\fin}{{\rm fin}}
\newcommand{\id}{{\rm id}}
\newcommand{\ie}{{\rm i}}
\newcommand{\GL}{{\rm GL}}
\newcommand{\gl}{{{\mathfrak g}{\mathfrak l}}}
\newcommand{\Hom}{{\rm Hom}\,}
\newcommand{\Img}{{\rm Im}\,}
\newcommand{\Ker}{{\rm Ker}\,}
\newcommand{\lf}{{\rm l}}
\newcommand{\Lie}{\text{\bf L}}
\newcommand{\m}{\textbf{m}}
\newcommand{\OO}{{\rm O}}
\newcommand{\pr}{{\rm pr}}
\newcommand{\Ran}{{\rm Ran}\,}
\renewcommand{\Re}{{\rm Re}\,}
\newcommand{\rb}{\text{\bf r}}
\newcommand{\Sp}{{\rm Sp}}
\renewcommand{\sp}{{{\mathfrak s}{\mathfrak p}}}
\newcommand{\spa}{{\rm span}\,}
\newcommand{\Tr}{{\rm Tr}\,}

\newcommand{\G}{{\rm G}}
\newcommand{\U}{{\rm U}}

\newcommand{\Ac}{{\mathcal A}}
\newcommand{\Bc}{{\mathcal B}}
\newcommand{\Cc}{{\mathcal C}}
\newcommand{\Ec}{{\mathcal E}}
\newcommand{\Dc}{{\mathcal D}}
\newcommand{\Hc}{{\mathcal H}}
\newcommand{\Jc}{{\mathcal J}}
\renewcommand{\Mc}{{\mathcal M}}
\newcommand{\Nc}{{\mathcal N}}
\newcommand{\Oc}{{\mathcal O}}
\newcommand{\Pc}{{\mathcal P}}
\newcommand{\Tc}{{\mathcal T}}
\newcommand{\Vc}{{\mathcal V}}
\newcommand{\Uc}{{\mathcal U}}

\newcommand{\Bg}{{\mathfrak B}}
\newcommand{\Fg}{{\mathfrak F}}
\newcommand{\Gg}{{\mathfrak G}}
\newcommand{\Ig}{{\mathfrak I}}
\newcommand{\Jg}{{\mathfrak J}}
\newcommand{\Lg}{{\mathfrak L}}
\newcommand{\Pg}{{\mathfrak P}}
\newcommand{\Sg}{{\mathfrak S}}
\newcommand{\Xg}{{\mathfrak X}}
\newcommand{\Yg}{{\mathfrak Y}}
\newcommand{\Zg}{{\mathfrak Z}}

\newcommand{\ag}{{\mathfrak a}}
\newcommand{\bg}{{\mathfrak b}}
\newcommand{\dg}{{\mathfrak d}}
\renewcommand{\gg}{{\mathfrak g}}
\newcommand{\hg}{{\mathfrak h}}
\newcommand{\kg}{{\mathfrak k}}
\newcommand{\mg}{{\mathfrak m}}
\newcommand{\n}{{\mathfrak n}}
\newcommand{\og}{{\mathfrak o}}
\newcommand{\pg}{{\mathfrak p}}
\newcommand{\sg}{{\mathfrak s}}
\newcommand{\tg}{{\mathfrak t}}
\newcommand{\ug}{{\mathfrak u}}

\makeatletter
\title[Iwasawa decompositions of some infinite-dimensional Lie groups]
{Iwasawa decompositions of some infinite-dimensional Lie groups}
\author{Daniel Belti\c t\u a}
\address{Institute of Mathematics ``Simion Stoilow'' 
of the Romanian Academy, 
RO-014700 Bucharest, Romania}
\email{Daniel.Beltita@imar.ro}
\date{January 15, 2007}
\makeatother

\begin{abstract} 
We set up an abstract framework 
that allows the investigation of Iwasawa decompositions 
for involutive infinite-dimensional Lie groups modeled on Banach spaces.
As an application, we construct Iwasawa decompositions for 
classical real or complex Banach-Lie groups 
associated with the Schatten ideals ${\mathfrak S}_p({\mathcal H})$ on 
a complex separable Hilbert space ${\mathcal H}$  
if $1<p<\infty$. 
\end{abstract}

\maketitle

\tableofcontents

\section{Introduction}\label{sectionIntr}

Our aim in this paper is to set up an abstract framework 
that allows the investigation of Iwasawa decompositions 
for involutive infinite-dimensional Lie groups modeled on Banach spaces. 
In particular we address an old conjecture on the existence 
of such decompositions for the classical Banach-Lie groups of operators 
associated with the Schatten operator ideals on Hilbert spaces 
(see subsection~8.4 of Section~II.8 in \cite{dlH72} 
and Section~\ref{sectionClass} below), 
and we show that the corresponding question can be answered 
in the affirmative in  many cases, 
even in the case of the covering groups. 

The Iwasawa decompositions of 
finite-dimensional reductive Lie groups 
(see e.g., \cite{Iwa49}, \cite{Hel01}, and \cite{Kna96}) 
play a crucial role in areas like 
differential geometry and  representation theory. 
So far as differential geometry is concerned, 
there exists a recent growth of interest in group decompositions 
and their implications in various geometric problems 
in infinite dimensions---see for instance 
\cite{Tum05}, \cite{Tum06}, and the references therein. 
From this point of view it is natural to try to understand 
the infinite-dimensional versions of Iwasawa decompositions as well; 
this problem was  already addressed in the case of loop groups 
in \cite{Kel99} and \cite{BD01}. 
As regards the representation theory,  it is well known that 
decompositions of this kind are 
particularly important for instance 
in the construction of principal series representations. 
And there has been a continuous endeavor 
to extend the ideas of representation theory 
to the setting of infinite-dimensional Lie groups. 
Some references related in spirit to the present paper are 
\cite{Se57}, \cite{Ki73}, \cite{SV75}, \cite{Ol78}
\cite{Boy80}, \cite{Ca85}, 
\cite{Ol88}, \cite{Pic90}, \cite{Boy93}, \cite{Nee98}, \cite{NO98}, 
\cite{NRW01}, \cite{DPW02}, \cite{Nee04}, \cite{Gru05}, 
\cite{Wo04}, \cite{BR06}, 
however this list is very far from being complete. 
In this connection we wish to highlight the paper \cite{Wo04} 
devoted to an investigation of direct limits of 
(Iwasawa decompositions and) principal series representations of 
reductive Lie groups. 
The results of the present paper can be thought of as belonging to 
the same line of investigation, 
inasmuch as the construction of Iwasawa decompositions 
should be  the very first step toward the construction of 
principal series representations for the classical Banach-Lie groups 
and their covering groups. 

Another source of interest in obtaining Iwasawa decompositions 
for infinite-dimensional versions of reductive Lie groups 
is related to the place held 
by reductive structures in the geometry of 
many infinite-dimensional manifolds---see for instance 
\cite{CG99} and \cite{Nee02c}. 
We refer also to the recent survey \cite{Ga06} 
which skilfully highlights the special relationship 
between the reductive structures and the idea of amenability. 
That relationship also plays an important role in the abstract framework 
constructed in Section~\ref{sectionAbs} of the present paper. 
It is noteworthy that reductive structures 
with a Lie theoretic flavor constitute the background of 
the papers \cite{Neu99} and \cite{Neu02} as well, 
concerning convexity theorems 
(cf.~\cite{Ko73}, \cite{LR91}, \cite{BFR93})
in an infinite-dimensional setting. 

The methods we use to set up the aforementioned abstract framework are 
 largely inspired by 
the interaction between the theory of Lie algebras and 
the local spectral theory of bounded operators (see \cite{BS01}). 
These methods turn out to be particularly effective in order 
to identify what the third component of 
an infinite-dimensional Iwasawa decomposition should be, 
that is, 
the infinite-dimensional versions of nilpotent Lie groups and algebras. 

On the other hand, the applications we make to 
the classical Banach-Lie groups are naturally related to 
the theory of triangular integrals 
(\cite{GK70}, \cite{EL72}, \cite{Erd78}, \cite{Ara78}, \cite{Dav88}) 
and to the theory of factorization of Hilbert space operators 
along nests of subspaces 
(\cite{Arv67}, \cite{GK70}, \cite{Erd72}, \cite{Lar85}, \cite{Pow86}, 
\cite{Pit88}, \cite{MSS88}). 
For the reader's convenience 
we recorded in Appendix~A some auxiliary facts on operator ideals,  
in particular
the factorization results suitable for our purposes. 

\addtocontents{toc}{\SkipTocEntry}
\subsection*{Rough decompositions}

Here we show a sample of pathological phenomenon 
one has to avoid in order to obtain \textit{smooth} Iwasawa decompositions 
for infinite-dimensional Lie groups.  

\begin{proposition}\label{rough}
Let $\Hc$ be a complex separable Hilbert space with an orthonormal basis 
$\{\xi_j\}_{0\le j<\omega}$, 
where $\omega\in{\mathbb N}\cup\{\aleph_0\}$. 
Then consider the Banach-Lie group 
$G:=\GL(\Hc)$ consisting of all invertible bounded linear operators 
on $\Hc$, and its subgroups 
$$\begin{aligned}
K:=&\{k\in G\mid k^*k=\1\}, \\
A:=&\{a\in G\mid
a\xi_j\in{\mathbb R}_{+}^*\xi_j\text{ whenever }0\le j<\omega\}, 
 \text{ and }\\
N:=&\{n\in G\mid
n\xi_j\in\xi_j+\spa\{\xi_l\mid l<j\}\text{ whenever }0\le j<\omega\}\}.
\end{aligned}$$
In addition, consider the Banach-Lie algebra $\gg=\Bc(\Hc)$, 
with its closed Lie subalgebras 
$$\begin{aligned}
\kg:=&\{X\in\gg\mid X^*=-X\}, \\
\ag:=&\{Y\in \gg\mid
Y\xi_j\in{\mathbb R}\xi_j\text{ whenever }0\le j<\omega\}, 
 \text{ and }\\
\n:=&\{Z\in\gg\mid
Z\xi_j\in\spa\{\xi_l\mid l<j\}\text{ whenever }0\le j<\omega\}\}.
\end{aligned}$$
Then the following assertions hold: 
\begin{itemize}
\item[$\bullet$]
$K$, $A$, and $N$ are Banach-Lie groups  
with the corresponding Lie algebras 
$\kg$, $\ag$, and $\n$, respectively, 
and the multiplication map 
$\m\colon K\times A\times N\to G$, $(k,a,n)\mapsto kan$,  
is smooth and bijective.
\item[$\bullet$] 
The mapping $\m$ 
is a diffeomorphism if and only if $\kg\dotplus\ag\dotplus\n=\gg$, 
and this equality holds if and only if 
the Hilbert space $\Hc$ is finite-dimensional. 
\end{itemize}
\end{proposition}

\begin{proof}
It is straightforward to prove that $\m$ is injective 
since $K\cap AN=\{\1\}$. 
To prove that $\m$ is surjective as well, 
we can use the unital Banach algebra 
$\Bg=
\{b\in\Bc(\Hc)\mid b\xi_j\in\spa\{\xi_l\mid 0\le l\le j\}
\text{ if }0\le j<\omega\}$. 
Denote by $\Bg^\times$ the group of invertible elements in $\Bg$. 
It was proved in \cite{Arv75} and \cite{Lar85} 
that for every $g\in\GL(\Hc)$ 
there exist $k\in K$ and $b\in\Bg^\times$ 
such that $g=kb$. 
It is easy to see that $\Bg^\times\subseteq KAN$, 
hence $g=kb\in KAN=\m(K\times A\times N)$. 
The fact that $K$, $A$, $N$ are Banach-Lie groups 
with the corresponding Lie algebras 
$\kg$, $\ag$, and $\n$, respectively, 
and 
follows for instance by Corollary~3.7 in \cite{Bel06}, 
and in addition the inclusion maps of $K$, $A$, and $N$ 
into $G$ are smooth. 
It then follows that 
the multiplication map 
$\m\colon K\times A\times N\to G$ 
is smooth as well. 

In order to prove the second assertion note that 
the tangent mapping 
$T_{(\1,\1,\1)}\m\colon\kg\times\ag\times\n\to\gg$ 
is given by $(X,Y,Z)\mapsto Z+Y+Z$, 
hence $\m$ is a local diffeomorphism if and only if 
$\kg\dotplus\ag\dotplus\n=\gg$. 
Since we have already seen that $\m$ is a bijective map, 
it follows that the latter direct sum decomposition actually holds 
if and only if $\m$ is a diffeomorphism. 
Next note that if $\dim\Hc<\infty$ then we get 
$\kg\dotplus\ag\dotplus\n=\gg$ by an elementary reasoning 
(or by the local Iwasawa decomposition for 
the complex finite-dimensional reductive Lie algebra $\gg$; 
see e.g., \cite{Kna96}). 

Thus, to complete the proof, 
it will be enough to show that if the Hilbert space $\Hc$ 
is infinite-dimensional, 
then $\gg\setminus(\kg+\ag+\n)\ne\emptyset$. 
In fact, for $j,l\in{\mathbb N}$ denote $h_{jl}=0$ if $j=l$ 
and $h_{jl}=1/(j-l)$ if $j\ne l$. 
Then there exists $W\in\Bc(\Hc)$ 
whose matrix with respect to the orthonormal basis 
$\{\xi_j\}_{j\in{\mathbb N}}$ is $(h_{jl})_{j,l\in{\mathbb N}}$,  
and in addition $W^*=-W$ and there exist no operators 
$Z_1,Z_2\in\n$ with $W=Z_1-Z_2^*$ 
(see Example~4.1 in~\cite{Dav88}). 
Now it is easy to see that $\ie W\not\in \kg+\ag+\n$. 
In fact, if $\ie W=X+Y+Z$ with $X\in\kg$, $Y\in\ag$, and $Z\in\n$, 
then $\ie W=(\ie W)^*=X^*+Y^*+Z^*=-X+Y+Z^*$. 
Thence $2\ie W=2Y+Z+Z^*$. 
Since the matrix of $W$ has only zeros on the diagonal, 
we get $Y=0$, whence $2\ie W=Z+Z^*$.  
Then $W=(-(\ie/2)Z)-(-(\ie/2)Z)^*$, 
and this contradicts one of the above mentioned properties 
of $W$. 
Thus $\ie W\in\gg\setminus(\kg+\ag+\n)$, and this completes the proof. 
\end{proof}

\section{Iwasawa decompositions for 
involutive Banach-Lie groups}\label{sectionAbs}

In this section we sketch an abstract framework 
that allows to investigate Iwasawa decompositions 
for involutive infinite-dimensional Lie groups modeled on Banach spaces. 
We will apply these abstract statements in 
Sections \ref{sectionA}, \ref{sectionB}, and \ref{sectionC} 
in the case of the classical Banach-Lie groups 
associated with norm ideals. 
The central idea of this abstract approach 
is that the local Iwasawa decompositions 
can be constructed out of certain special elements 
of Lie algebras, which we call Iwasawa regular elements 
(Definition~\ref{iw_abstr1}). 

\addtocontents{toc}{\SkipTocEntry}
\subsection*{Preliminaries on local spectral theory}

Throughout the paper we let  
${\mathbb R}$, ${\mathbb C}$, and ${\mathbb H}$ 
stand for 
the fields of the real, complex, and quaternionic numbers, 
respectively. 

For a real or complex Banach space ${\mathfrak X}$ we denote 
either by $\id_{\mathfrak X}$ or simply by $\1$ 
the identity map of ${\mathfrak X}$, 
by ${\mathfrak X}^*$ the topological dual of ${\mathfrak X}$,   
by ${\mathcal B}({\mathfrak X})$ the algebra of 
all bounded linear operators on ${\mathfrak X}$ and, 
when ${\mathfrak X}$ is a complex Banach space,  
we denote by $\sigma(D)$ the spectrum of $D$ 
whenever $D\in{\mathcal B}({\mathfrak X})$ . 
In this case, for every $x\in{\mathfrak X}$ 
we denote by $\sigma_D(x)$ the \textit{local spectrum} of $x$ 
with respect to $D$. 
We recall that $\sigma_D(x)$ is a closed subset of $\sigma(D)$ and  
$w\in{\mathbb C}\setminus\sigma_D(x)$ 
if and only if 
there exists an open neighborhood $W$ of $w$ and a holomorphic function 
$\xi\colon W\to{\mathfrak X}$ such that 
$(z\id_{\mathfrak X}-D)\xi(z)=x$ for every $z\in W$.  
If $F\subseteq{\mathbb C}$ we denote 
$${\mathfrak X}_D(F)=\{x\in{\mathfrak X}\mid\sigma_D(x)\subseteq F\}.$$ 
We note that, in the case when ${\mathfrak X}$ has {\it finite dimension} $m$, 
we have
$${\mathfrak X}_D(F)=
\bigoplus_{\lambda\in F\cap\sigma(D)}\Ker((D-\lambda\id_{\mathfrak X})^m)
\quad\text{ for every }F\subseteq{\mathbb C},$$ 
while if ${\mathfrak X}$ is a Hilbert space 
and $D$ is a normal operator with the spectral measure $E_D(\cdot)$,  
then 
$${\mathfrak X}_D(F)=\Ran E_D(F)
\quad\text{ whenever $F$ is a closed subset of }{\mathbb C}.$$ 
See \S 12 in \cite{BS01} for a review of the few facts 
needed from the local spectral theory. 
(More bibliographical details can be found in 
the Notes to Chapter~I in \cite{BS01}.)

\addtocontents{toc}{\SkipTocEntry}
\subsection*{Projections on kernels of skew-Hermitian operators}

\begin{notation}\label{D1}
\normalfont
The following notation will be used throughout the paper:
\begin{itemize}
\item[$\bullet$] 
For every complex Banach space $\Xg$ we denote 
$\ell^\infty_{\Xg}({\mathbb R}):=
\{f\colon{\mathbb R}\to{\Xg}\mid \Vert f\Vert_{\infty} 
:=\sup\limits_{\mathbb R}\Vert f(\cdot)\Vert<\infty\}$, 
which is in turn a complex Banach space. 
\item[$\bullet$] We pick a state 
$\mu\colon\ell^\infty_{\mathbb C}({\mathbb R})\to{\mathbb C}$ 
of the commutative unital $C^*$-algebra 
$\ell^\infty_{\mathbb C}({\mathbb R})$ satisfying 
the translation invariance condition 
$$(\forall f\in\ell^\infty_{\mathbb C}({\mathbb R}))\,
(\forall t\in{\mathbb R})\quad 
\mu(f(\cdot))=\mu(f(\cdot+t)) $$
and the symmetry condition 
$$(\forall f\in\ell^\infty_{\mathbb C}({\mathbb R}))\quad
\mu(f(\cdot))=\mu(f(-\cdot)). $$
(See Problem~7 in Chapter~2 of \cite{Pat88}.) 
For every $f\in\ell^\infty_{\mathbb C}({\mathbb R})$ 
we denote 
$\mu(f(\cdot))=\int\limits_{\mathbb R}f(t)\de\mu(t)$.  
\item[$\bullet$] For every complex Banach space $\Xg$ 
and every $f\in\ell^\infty_{\Xg}({\mathbb R})$ we define 
$\mu(f)=\int\limits_{\mathbb R}f(t)\de\mu(t)\in\Xg^{**}$ 
by the formula 
$$(\mu(f))(\varphi)=\int\limits_{\mathbb R}\varphi(f(t))\de\mu(t)$$ 
whenever $\varphi\in\Xg^{*}$ (see \cite{BP05}). 
\end{itemize}
\qed
\end{notation}

\begin{definition}\label{D2}
\normalfont
Let $\Xg_0$ be a real Banach space 
with the complexification $\Xg=(\Xg_0)^{\mathbb C}=\Xg_0+\ie\Xg_0$, 
which is a complex Banach space with the norm given by 
$\|y_1+\ie y_2\|
:=\sup\limits_{t\in[0,2\pi]}\|(\cos t)y_1+(\sin t)y_2\|$ 
for all $y_1,y_2\in\Xg_0$. 
Assume that $A\colon\Xg_0\to\Xg_0$ is a bounded linear operator 
such that $\sup\limits_{t\in\mathbb R}\Vert\exp(tA)\Vert<\infty$, 
and denote by $A\colon\Xg\to\Xg$ the ${\mathbb C}$-linear extension 
of $A\colon\Xg_0\to\Xg_0$. 
In this case we define 
$$\Xg_{A}^{+}:=
\Xg_A(\ie[0,\infty))=\Bigl\{y\in\Xg\mid
\limsup\limits_{t\to\infty}\frac{1}{t}\log\Vert(\exp(\ie tA))y\Vert\le0
\Bigr\}$$
(see also Remark~1.3 in \cite{Be01}). 
Now assume that $\Xg_0$ is a \textit{reflexive} Banach space. 
Then $\Xg$ will be a reflexive complex Banach space, 
and there exists a bounded linear operator 
$\Dc_{\Xg,A}\colon\Xg\to\Xg$ 
such that 
$$(\forall y\in\Xg)\quad
\Dc_{\Xg,A}y=\int\limits_{\mathbb R}(\exp(tA))y\,\de\mu(t).$$
It is easy to see that $\Dc_{\Xg,A}\Xg_0\subseteq\Xg_0$, 
and we shall define 
$\Dc_{\Xg_0,A}:=\Dc_{\Xg,A}|_{\Xg_0}\colon\Xg_0\to\Xg_0$. 
 
In addition,
we define
$\Xg_{A}^{0,+}:=\Xg_{A}^{+}\cap\Ker\Dc_{\Xg,A}$.
\qed
\end{definition}

\begin{remark}\label{D3}
\normalfont
In the setting of Definition~\ref{D2} we have 
$(\Dc_{\Xg,A})^2=\Dc_{\Xg,A}$, $\Ran\Dc_{\Xg,A}=\Ker A\subseteq\Xg$, 
and $\Ran\Dc_{\Xg_0,A}=(\Ker A)\cap\Xg_0$ 
(see \cite{BP05}). 
\qed
\end{remark}

\addtocontents{toc}{\SkipTocEntry}
\subsection*{Elliptic involutive Banach-Lie algebras and 
abstract Iwasawa decompositions}

\begin{definition}\label{ell1}
\normalfont
Let $\gg_0$ be an involutive real or complex Banach-Lie algebra, that is, 
$\gg_0$ is equipped with a continuous linear mapping $X\mapsto X^*$ 
such that $(X^*)^*=X$ and $[X,Y]^*=-[X^*,Y^*]$ 
whenever $X,Y\in\gg_0$. 
If $\gg_0$ is a complex Banach-Lie algebra, then we assume in addition 
that $(\ie X)^*=-\ie X^*$ for all $X\in\gg_0$. 

We say that $\gg_0$ is an \textit{elliptic} involutive Banach-Lie algebra 
if $\Vert\exp(t\cdot\ad_{\gg_0}X)\Vert\le1$ whenever 
$t\in{\mathbb R}$ and $X=-X^*\in\gg_0$. 
\qed
\end{definition}

\begin{remark}\label{ell2}
\normalfont
In the special case of the canonically involutive real Banach-Lie algebras 
(that is, 
 $X^*=-X$ for all $X\in\gg_0$) the above Definition~\ref{ell1} 
coincides with Definition~IV.3 in \cite{Nee02c}
(or Definition~8.24 in \cite{Bel06}). 
\qed
\end{remark}

\begin{definition}\label{iw_abstr1}
\normalfont
Let $\gg_0$ be an elliptic real Banach-Lie algebra 
with the complexification $\gg$, 
and denote $\pg_0:=\{X\in\gg_0\mid X^*=X\}$. 
An \textit{Iwasawa decomposition of $\gg_0$} 
is a direct sum decomposition 
\begin{equation}\label{iw_local}
\gg_0=\kg_0\dotplus\ag_0\dotplus\n_0
\end{equation}
satisfying the following conditions:
\begin{itemize}
\item[{\rm(j)}] We have 
$\kg_0=\{X\in\gg_0\mid X^*=-X\}$.  
\item[{\rm(jj)}] The term $\ag_0$ is a linear subspace of 
$\pg_0$ such that $[\ag_0,\ag_0]=\{0\}$. 
\item[{\rm(jjj)}] There exists $X_0\in\ag_0$ such that 
$\ag_0=\ag_{X_0}$ and $\n_0=\n_{X_0}$, 
where 
$\ag_{X_0}=\pg_0\cap\Ker(\ad\,X_0)$ and 
$\n_{X_0}=\gg_0\cap\gg_{\ad(-\ie X_0)}^{0,+}$. 
\end{itemize}
In this case we say that $X_0$ is an \textit{Iwasawa regular element} 
of $\gg_0$ and \eqref{iw_local} is the 
\textit{Iwasawa decomposition of $\gg_0$ associated with $X_0$}. 
In the case when all of the conditions \eqref{iw_local} 
and (j)--(jjj) are satisfied perhaps except for (jj), 
we say that $X_0$ is an \textit{Iwasawa quasi-regular element} 
(and \eqref{iw_local} is still called the 
Iwasawa decomposition of $\gg_0$ associated with $X_0$). 

Now let us assume that  
$G$ is a connected Banach-Lie group with $\Lie(G)=\gg_0$ 
and $K$, $A$, and $N$ are the connected Banach-Lie groups 
which are subgroups of $G$ and correspond to the Lie algebras 
$\kg_0$, $\ag_0$, and $\n_0$, respectively. 
If the mapping 
\begin{equation}\label{iw_global}
\m\colon K\times A\times N\to G,\quad (k,a,n)\mapsto kan 
\end{equation}
is a diffeomorphism, 
then we say that this mapping is 
the \textit{global Iwasawa decomposition of $G$} corresponding 
to~\eqref{iw_local}. 
\qed
\end{definition}

\begin{remark}\label{iw_abstr1.5}
\normalfont
In the setting of Definition~\ref{iw_abstr1}, 
if $X_0$ is an Iwasawa regular element then 
it is easy to see that $\ag_{X_0}$ is a \textit{maximal} linear subspace of 
$\pg_0$ such that $[\ag_{X_0},\ag_{X_0}]=\{0\}$, 
while $\n_{X_0}$ is a closed subalgebra of $\gg_0$ 
(see also \cite{Be01}). 
\qed
\end{remark}

\begin{remark}\label{iw_abstr2}
\normalfont
In the setting of Definition~\ref{iw_abstr1}, 
it is easy to see that all of the groups $K$, $A$, and $N$ 
are Banach-Lie subgroups of $G$. 
(See \cite{Up85} or \cite{Bel06} for details on the latter notion.) 
\qed
\end{remark}

\begin{proposition}\label{iw_abstr3}
In the setting of {\rm Definition~\ref{iw_abstr1}}, 
let us assume that the Banach-Lie algebra $\gg_0$ 
is actually an elliptic involutive complex Banach-Lie algebra. 
Then for every $X\in\pg_0$ and every closed subset $F$ of ${\mathbb R}$ 
we have 
$\gg_0\cap\gg_{\ad_{\gg}X}(F)
=(\gg_0)_{\ad_{\gg_0}X}(F)$. 
\end{proposition}

\begin{proof}
Since $\gg_0$ is an elliptic Banach-Lie algebra, 
it follows that $\ad_{\gg_0}X\colon\gg_0\to\gg_0$ is 
a Hermitian operator (see Definition~5.23 in \cite{Bel06}). 
If $\gg$ stands for the complexification of $\gg_0$, 
then  $\ad_{\gg}X\colon\gg\to\gg$ is Hermitian as well. 
In particular, there exist quasimultiplicative maps 
$\Psi_{\ad_{\gg_0}X}\colon\Ci({\mathbb R})\to\Bc(\gg_0)$ and 
$\Psi_{\ad_{\gg}X}\colon\Ci({\mathbb R})\to\Bc(\gg)$ such that 
$\Psi_{\ad_{\gg_0}X}(\id_{\mathbb R})=\ad_{\gg_0}X$ and 
$\Psi_{\ad_{\gg}X}(\id_{\mathbb R})=\ad_{\gg}X$, respectively. 
The maps $\Psi_{\ad_{\gg_0}X}$ and $\Psi_{\ad_{\gg}X}$ 
can be constructed by the Weyl functional calculus 
as in Example~5.25 in \cite{Bel06}. 

Now let $\iota\colon\gg_0\hookrightarrow\gg$ be the inclusion map. 
Then $\iota(X)=X$, hence Remark~5.19 in \cite{Bel06} shows that  
for every closed subset $F$ of ${\mathbb R}$ 
we have $\iota((\gg_0)_{\ad_{\gg_0}X}(F))\subseteq\gg_{\ad_{\gg}X}(F)$, 
whence 
$(\gg_0)_{\ad_{\gg_0}X}(F)\subseteq\gg_0\cap\gg_{\ad_{\gg}X}(F)$.

To prove the converse inclusion, denote by $\kappa\colon\gg\to\gg$, 
the conjugation on $\gg$ whose fixed point set 
is $\gg_0$, and then define $\pi\colon\gg\to\gg_0$, $\pi(Z)=(Z+\kappa(Z))/2$. 
Then $\pi(X)=X$, 
whence $\pi\circ(\ad_{\gg}X)=(\ad_{\gg_0}X)\circ\pi$. 
Now Remark~5.19 in \cite{Bel06} again shows that 
for every closed subset $F$ of ${\mathbb R}$ 
we have $\pi(\gg_{\ad\,X}(F))\subseteq\gg_{\ad_{\gg_0}X}(F)$. 
Thence $\gg_0\cap\gg_{\ad\,X}(F)\subseteq\gg_{\ad_{\gg_0}X}(F)$, 
and we are done. 
\end{proof}

\begin{remark}\label{iw_abstr4}
\normalfont
In the special case when $\dim\gg_0<\infty$ and  $F$ 
is a certain finite subset of ${\mathbb R}_{+}$, 
the conclusion of our Proposition~\ref{iw_abstr3} 
was obtained in Chapter~VI, \S 6 of \cite{Hel01} 
by using the structure theory of 
finite-dimensional complex semisimple Lie algebras. 
\qed
\end{remark}

\begin{proposition}\label{D4}
Let $\widetilde{\gg}$ be an elliptic complex Banach-Lie algebra 
whose underlying Banach space is reflexive, 
and pick $X_0=X_0^*\in\widetilde{\gg}$.  
Assume that we have a bounded linear operator 
$\widetilde{\Tc}\in\Bc(\widetilde{\gg})$ 
such that 
\begin{eqnarray}
\widetilde{\Tc}^2&=&\widetilde{\Tc},\nonumber \\
\label{D**}
\Ran\widetilde{\Tc}&=&\widetilde{\gg}_{\ad\,X_0}({\mathbb R}_{+}),\\
\label{D***}
\Ran(\1-\widetilde{\Tc})&\subseteq& 
 \widetilde{\gg}_{\ad\,X_0}(-{\mathbb R}_{+}).
\end{eqnarray}

Then $X_0$ is an Iwasawa quasi-regular element of $\widetilde{\gg}$. 
Let 
\begin{equation}\label{D*}
\widetilde{\gg}=\widetilde{\kg}\dotplus\widetilde{\ag}\dotplus\widetilde{\n}
\end{equation}
be the Iwasawa decomposition of $\widetilde{\gg}$ associated with $X_0$,  
and for 
$\widetilde{\sg}\in\{\widetilde{\kg},\widetilde{\ag},\widetilde{\n}\}$, 
denote by $p_{\widetilde{\sg}}\colon\widetilde{\gg}\to\widetilde{\sg}$ 
the linear projections corresponding to 
the direct sum decomposition~\eqref{D*}. 
Then for all $X\in\widetilde{\gg}$ we have 
\begin{equation*}
\begin{aligned}
p_{\widetilde{\kg}}(X)
  &=(\1-\widetilde{\Tc})X-((\1-\widetilde{\Tc})X)^*
    +\frac{1}{2}(\Dc_{\widetilde{\gg},\ad(-\ie X_0)}X
     -(\Dc_{\widetilde{\gg},\ad(-\ie X_0)}X)^*), \\
p_{\widetilde{\ag}}(X)
  &=\frac{1}{2}(\Dc_{\widetilde{\gg},\ad(-\ie X_0)}X
     +(\Dc_{\widetilde{\gg},\ad(-\ie X_0)}X)^*), \text{ and }\\
p_{\widetilde{\n}}(X)
  &=(\widetilde{\Tc}-\Dc_{\widetilde{\gg},\ad(-\ie X_0)})X
    +((\1-\widetilde{\Tc})X)^*.    
\end{aligned}
\end{equation*}
\end{proposition}

\begin{proof}
To begin with, recall that 
\begin{eqnarray}
\label{D(1)} 
\widetilde{\kg}
 &=&\{X\in\widetilde{\gg}\mid X^*=-X\}, \\
\label{D(2)}
\widetilde{\ag}
 &=&\{X\in\widetilde{\gg}\mid X^*=X\text{ and }[X_0,X]=0\}, \text{ and }\\
\label{D(3)}
\widetilde{\n}
 &=&\widetilde{\gg}_{\ad(-\ie X_0)}^{0,+}
  =\widetilde{\gg}_{\ad\,X_0}({\mathbb R}_{+})
    \cap\Ker\Dc_{\widetilde{\gg},\ad(-\ie X_0)},
\end{eqnarray}
where the latter two equalities follow by 
Proposition~\ref{iw_abstr3}, Definition~\ref{D2}, 
and Definition~\ref{iw_abstr1}. 
It is straightforward to check that 
$\widetilde{\kg}\cap(\widetilde{\ag}+\widetilde{\n})
=\widetilde{\ag}\cap\widetilde{\n}=\{0\}$, 
whence 
$$\widetilde{\kg}\cap(\widetilde{\ag}+\widetilde{\n})
=\widetilde{\ag}\cap(\widetilde{\kg}+\widetilde{\n})
=\widetilde{\n}\cap(\widetilde{\kg}+\widetilde{\ag})
=\{0\},$$
and it remains to prove that 
$\widetilde{\kg}+\widetilde{\ag}+\widetilde{\n}=\widetilde{\gg}$. 

For this purpose, first note that $\Dc_{\widetilde{\gg},\ad(-\ie X_0)}$ 
and $\widetilde{\Tc}$ are idempotent operators on $\widetilde{\gg}$ 
satisfying 
$$\Ran\Dc_{\widetilde{\gg},\ad(-\ie X_0)}\subseteq\Ran\widetilde{\Tc},$$ 
hence 
\begin{equation}\label{D(4)} 
\widetilde{\Tc}\Dc_{\widetilde{\gg},\ad(-\ie X_0)}
=\Dc_{\widetilde{\gg},\ad(-\ie X_0)}\widetilde{\Tc}
=\Dc_{\widetilde{\gg},\ad(-\ie X_0)}. 
\end{equation}
Now let $X\in\widetilde{\gg}$ arbitrary and denote by 
$X_{\kg}$, $X_{\ag}$, and $X_{\n}$ 
the right-hand sides of the wished-for formulas 
for $p_{\widetilde{\kg}}(X)$, $p_{\widetilde{\ag}}(X)$, 
and $p_{\widetilde{\n}}(X)$, respectively. 
Thus we have to prove that 
$p_{\widetilde{\sg}}(X)=X_{\sg}$ 
for $\sg\in\{\kg,\ag,\n\}$. 
Moreover, it is clear that $X=X_{\kg}+X_{\ag}+X_{\n}$, 
so it will be enough to check that 
$X_{\kg}\in\widetilde{\kg}$, $X_{\ag}\in\widetilde{\ag}$, 
and $X_{\n}\in\widetilde{\n}$. 

It follows at once by \eqref{D(1)} that $X_{\kg}\in\widetilde{\kg}$. 
To see that $X_{\ag}\in\widetilde{\ag}$, 
first note that $X_{\ag}^*=X_{\ag}$. 
On the other hand, we have 
$[X_0,\Dc_{\widetilde{\gg},\ad(-\ie X_0)}X]=0$ 
(see Remark~\ref{D3}). 
Since $X_0=X_0^*$, 
it then follows that $[X_0,(\Dc_{\widetilde{\gg},\ad(-\ie X_0)}X)^*]=0$, 
whence $[X_0,X_{\ag}]=0$. 
Thus $X_{\ag}\in\widetilde{\ag}$. 

It remains to show that $X_{\n}\in\widetilde{\n}$. 
To this end, first note that 
$$\Dc_{\widetilde{\gg},\ad(-\ie X_0)}
(\widetilde{\Tc}-\Dc_{\widetilde{\gg},\ad(-\ie X_0)})X=0
\quad\text{and}\quad  
\Dc_{\widetilde{\gg},\ad(-\ie X_0)}(((\1-\widetilde{\Tc})X)^*)
=(\Dc_{\widetilde{\gg},\ad(-\ie X_0)}(\1-\widetilde{\Tc})X)^*=0 $$
by \eqref{D(4)} and the fact that 
$(\Dc_{\widetilde{\gg},\ad(-\ie X_0)}Y)^2
=\Dc_{\widetilde{\gg},\ad(-\ie X_0)}Y$. 
(We also used the fact that 
$\Dc_{\widetilde{\gg},\ad(-\ie X_0)}(Y^*)
=(\Dc_{\widetilde{\gg},\ad(-\ie X_0)}Y)^*$ 
whenever $Y\in\widetilde{\gg}$, 
which is a 
consequence of Definition~\ref{D2} 
since  
$X_0^*=X_0$.) 
It then follows that $\Dc_{\widetilde{\gg},\ad(-\ie X_0)}(X_{\n})=0$. 
Thus, according to~\eqref{D(3)}, we still have to prove that 
$X_{\n}\in\widetilde{\gg}_{\ad(-\ie X_0)}({\mathbb R}_{+})$. 
For this purpose we are going to show that both terms in 
the expression of $X_{\n}$ belong to 
$\widetilde{\gg}_{\ad\,X_0}({\mathbb R}_{+})$. 
In fact, by \eqref{D(4)}~and~\eqref{D**} we get 
$(\widetilde{\Tc}-\Dc_{\widetilde{\gg},\ad(-\ie X_0)})X
=\widetilde{\Tc}(\1-\Dc_{\widetilde{\gg},\ad(-\ie X_0)})X
\in\Ran\widetilde{\Tc}\subseteq\widetilde{\gg}_{\ad\,X_0}({\mathbb R}_{+})$. 
On the other hand, 
the mapping 
$\kappa\colon\widetilde{\gg}\to\widetilde{\gg}$, $Y\mapsto Y^*$ 
has the property $\theta\circ\ad\,X_0=-\ad\,X_0\circ\theta$ 
since $X_0^*=X_0$. 
Then Proposition~5.22 in \cite{Bel06} 
shows that 
$\kappa(\widetilde{\gg}_{\ad\,X_0}(-{\mathbb R}_{+}))
\subseteq\widetilde{\gg}_{\ad\,X_0}({\mathbb R}_{+})$, 
whence by \eqref{D***} we get 
$((\1-\widetilde{\Tc})X)^*=\kappa((\1-\widetilde{\Tc})X)
\in\widetilde{\gg}_{\ad\,X_0}({\mathbb R}_{+})$. 
Consequently 
$$X_{\n}=(\widetilde{\Tc}-\Dc_{\widetilde{\gg},\ad(-\ie X_0)})X
    +((\1-\widetilde{\Tc})X)^*
\in\widetilde{\gg}_{\ad\,X_0}({\mathbb R}_{+}),$$ 
and the proof is complete.  
\end{proof}

\begin{corollary}\label{D5}
Assume the setting of {\rm Proposition~\ref{D4}} and let 
$\gg$ be a closed involutive complex subalgebra of 
$\widetilde{\gg}$ such that 
$$X_0\in\gg\quad{and}\quad\widetilde{\Tc}(\gg)\subseteq\gg.$$
Then $X_0$ is an Iwasawa quasi-regular element of $\gg$ and 
the Iwasawa decomposition of $\gg$ associated with~$X_0$ 
is 
$\gg=(\widetilde{\kg}\cap\gg)\dotplus(\widetilde{\ag}\cap\gg)
\dotplus(\widetilde{\n}\cap\gg)$. 
\end{corollary}

\begin{proof}
Since $\widetilde{\Tc}(\gg)\subseteq\gg$ and 
$\Dc_{\widetilde{\gg},\ad(-\ie X_0)}(\gg)\subseteq\gg$, 
it follows by Proposition~\ref{D4} 
that $p_{\widetilde{\kg}}(\gg)\subseteq\gg$, 
$p_{\widetilde{\ag}}(\gg)\subseteq\gg$, 
and $p_{\widetilde{\n}}(\gg)\subseteq\gg$. 
It then follows by the direct sum decomposition~\eqref{D*} 
that $\gg=(\widetilde{\kg}\cap\gg)\dotplus(\widetilde{\ag}\cap\gg)
\dotplus(\widetilde{\n}\cap\gg)$. 
Thus, to conclude the proof, it remains to prove that 
$$\widetilde{\ag}\cap\gg=\{X\in\gg\mid X^*=X\text{ and }[X,X_0]=0\}
\quad\text{ and }\quad
\widetilde{\n}\cap\gg=\gg_{\ad(-\ie X_0)}^{0,+}$$
(see also Proposition~\ref{iw_abstr3}). 
The equality involving $\widetilde{\ag}\cap\gg$ is obvious. 
To prove the equality involving $\widetilde{\n}\cap\gg$, 
just note that by Proposition~\ref{iw_abstr3} we have 
$$\begin{aligned}
\widetilde{\n}\cap\gg 
 &=\widetilde{\gg}_{\ad(-\ie X_0)}^{0,+}\cap\gg 
  =\widetilde{\gg}_{\ad(-\ie X_0)}({\mathbb R}_{+})\cap
   (\Ker\Dc_{\widetilde{\gg},\ad(-\ie X_0)})\cap\gg \\
 &=\gg_{\ad(-\ie X_0)}({\mathbb R}_{+})\cap
   (\Ker\Dc_{\gg,\ad(-\ie X_0)}) 
 =\gg_{\ad(-\ie X_0)}^{0,+} 
\end{aligned}
$$
and this completes the proof. 
\end{proof}

\begin{corollary}\label{D6}
Assume the setting of {\rm Proposition~\ref{D4}} and let 
$\gg_0$ be a closed involutive real subalgebra of 
$\widetilde{\gg}$ such that 
$$X_0\in\gg_0,\quad\widetilde{\Tc}(\gg_0)\subseteq\gg_0, 
\quad\text{and}\quad\gg_0\cap\ie\gg_0=\{0\}.$$
Then $X_0$ is an Iwasawa quasi-regular element of $\gg$ and 
the Iwasawa decomposition of $\gg$ associated with~$X_0$ 
is 
$\gg_0=(\widetilde{\kg}\cap\gg_0)\dotplus(\widetilde{\ag}\cap\gg_0)
\dotplus(\widetilde{\n}\cap\gg_0)$. 
\end{corollary}

\begin{proof}
Denote $\gg:=\gg_0+\ie\gg_0$ ($\subseteq\widetilde{\gg}$). 
Since $\gg_0\cap\ie\gg_0=\{0\}$, 
it follows that $\gg$ is isomorphic 
to the complexification of $\gg_0$ 
(as a complex involutive Banach-Lie algebra). 

We have $X_0\in\gg_0\subseteq\gg$ and 
$\widetilde{\Tc}(\gg)=\widetilde{\Tc}(\gg_0+\ie\gg_0)\subseteq\gg_0+\ie\gg_0
=\gg$, 
hence Corollary~\ref{D5} shows that $X_0$ is Iwasawa regular in $\gg$ 
and the corresponding Iwasawa decomposition of $\gg$ 
is 
$\gg=(\widetilde{\kg}\cap\gg)\dotplus(\widetilde{\ag}\cap\gg)
\dotplus(\widetilde{\n}\cap\gg)$. 

In particular we get 
$\gg_{\ad(-\ie X_0)}^{0,+}=\widetilde{\n}\cap\gg$, 
whence 
$$\gg_0\cap\gg_{\ad(-\ie X_0)}^{0,+}=\widetilde{\n}\cap\gg_0.$$
On the other hand, it is obvious that 
$$\{X\in\gg_0\mid X^*=-X\}=\widetilde{\kg}\cap\gg_0
\quad\text{and}\quad 
\{X\in\gg_0\mid X^*=X\text{ and }[X_0,X]=0\}=\widetilde{\ag}\cap\gg_0, 
$$
hence the wished-for conclusion will follow as soon as we prove 
that 
$\gg_0=(\widetilde{\kg}\cap\gg_0)\dotplus(\widetilde{\ag}\cap\gg_0)
\dotplus(\widetilde{\n}\cap\gg_0)$. 
And this direct sum decomposition 
can be obtained just as in the proof of Corollary~\ref{D5}. 
Indeed, we have $\widetilde{\Tc}(\gg_0)\subseteq\gg_0$ 
and $\Dc_{\widetilde{\gg},\ad(-ie X_0)}(\gg_0)\subseteq\gg_0$. 
Then we can use Proposition~\ref{D4} to show that 
$p_{\widetilde{\kg}}(\gg_0)\subseteq\gg_0$, 
$p_{\widetilde{\ag}}(\gg_0)\subseteq\gg_0$, 
and $p_{\widetilde{\n}}(\gg_0)\subseteq\gg_0$. 
Since 
$\widetilde{\gg}=\widetilde{\kg}\dotplus\widetilde{\ag}\dotplus\widetilde{\n}$ 
by the hypothesis, 
it then follows that 
$\gg_0=(\widetilde{\kg}\cap\gg_0)\dotplus(\widetilde{\ag}\cap\gg_0)
\dotplus(\widetilde{\n}\cap\gg_0)$.
\end{proof}

\addtocontents{toc}{\SkipTocEntry}
\subsection*{Inductive limits of Iwasawa decompositions}

\begin{lemma}\label{top}
Let $\widetilde{\Psi}\colon\widetilde{S}_1\to\widetilde{S}_2$ 
be an open bijective mapping between two topological spaces. 
Assume that $S_j$ is a closed subset of $\widetilde{S}_j$ for $j=1,2$ 
such that $\widetilde{\Psi}(S_1)$ is a dense subset of $S_2$. 
Then $\widetilde{\Psi}(S_1)=S_2$.  
\end{lemma}

\begin{proof}
We have to prove that $S_2\subseteq\widetilde{\Psi}(S_1)$. 
The hypothesis that $\widetilde{\Psi}\colon\widetilde{S}_1\to\widetilde{S}_2$  
is an open bijection implies that its inverse 
$\widetilde{\Psi}^{-1}\colon\widetilde{S}_2\to\widetilde{S}_1$ 
is continuous. 
Then by using the other hypothesis, 
namely $\overline{\widetilde{\Psi}(S_1)}=S_2$,  
we get 
$ \widetilde{\Psi}^{-1}(S_2)
=\widetilde{\Psi}^{-1}(\overline{\widetilde{\Psi}(S_1)})
\subseteq\overline{\widetilde{\Psi}^{-1}(\widetilde{\Psi}(S_1))}
=\overline{S_1}=S_1$, 
which concludes the proof since $\widetilde{\Psi}$ is a bijection. 
\end{proof}

\begin{proposition}\label{lim_iw}
Let $\widetilde{G}$ be a Banach-Lie group and assume that 
$\widetilde{K}$, $\widetilde{A}$, and $\widetilde{N}$ are 
Banach-Lie subgroups of $\widetilde{G}$ such that 
the multiplication map 
$\widetilde{\m}\colon
\widetilde{K}\times\widetilde{A}\times\widetilde{N}\to\widetilde{G}$,
$(\widetilde{k},\widetilde{a},\widetilde{n})\mapsto
\widetilde{k}\widetilde{a}\widetilde{n}
$
is a diffeomorphism. 

Then let $G$, $K$, $A$, and $N$ be 
four connected Lie subgroups of $\widetilde{G}$ 
with $\Lie(G)=\gg$, $\Lie(K)=\kg$, $\Lie(A)=\ag$, and $\Lie(N)=\n$, 
and assume that $\gg$ is an elliptic real Banach-Lie algebra 
and $X_0\in\gg$ is an Iwasawa regular element 
such that the following conditions are satisfied: 
\begin{itemize}
\item[{\rm(j)}] 
We have $K\subseteq\widetilde{K}\cap G$, $A\subseteq\widetilde{A}\cap G$, 
$N\subseteq\widetilde{N}\cap G$, 
and $G$ is a Banach-Lie subgroup of $\widetilde{G}$.
\item[{\rm(jj)}] 
The Iwasawa decomposition of $\gg$ with respect to $X_0$ is 
$\gg=\kg\dotplus\ag\dotplus\n$.  
\item[{\rm(jjj)}] 
There exists a family $\{\gg_i\}_{i\in I}$ consisting 
of finite-dimensional reductive subalgebras of $\gg$ 
such that 
\begin{itemize}
\item[$\bullet$] 
there exists a bounded linear 
map $\Ec_i\colon\gg\to\gg$ such that $\Ec_i(X^*)=\Ec_i(X)^*$ for all $X\in\gg$, 
\item[$\bullet$]
$(\Ec_i)^2=\Ec_i$, $\Ran \Ec_i=\gg_i$, and $[\Ran(\1-\Ec_i),\gg_i]=\{0\}$, 
\item[$\bullet$]
$\Ec_i(X_0)$ is an Iwasawa regular element of $\gg_i$, and 
\item[$\bullet$] 
the connected subgroup $G_i$ of $G$ with $\Lie(G_i)=\gg_i$ 
is a closed subgroup and a finite-dimensional reductive Lie group
\end{itemize} 
for every $i\in I$, and $\overline{\bigcup\limits_{i\in I}\gg_i}=\gg$. 
\end{itemize}
Then the mapping 
$\m:=\widetilde{\m}|_{K\times A\times N}\colon K\times A\times N\to G$ 
is a diffeomorphism, $AN=NA$, 
and both $A$ and $N$ are simply connected. 
\end{proposition}

\begin{proof} 
Since $\widetilde{\m}$ is smooth and $G$ is a Banach-Lie subgroup of 
$\widetilde{G}$ 
by hypothesis~(j), 
it follows that $\m$ is a smooth mapping. 
Then condition~(jj) shows that the tangent map of $\m$ at any point 
of $K\times A\times N$ is an invertible continuous linear operator, 
hence $\m$ is a local diffeomorphism. 
On the other hand, $\m$ is injective since $\widetilde{\m}$ is so. 

It remains to prove that $\m$ is surjective. 
To this end let $i\in I$ and 
denote by $\gg_i=\kg_i\dotplus\ag_i\dotplus\n_i$ 
the Iwasawa decomposition of $\gg_i$ with respect to 
$\Ec_i(X_0)$. 
Also let $G_i=K_iA_iN_i$ be the corresponding global Iwasawa decomposition 
of the finite-dimensional reductive Lie group $G_i$ 
(see \cite{Kna96}). 
We have 
$$\kg_i=\{X\in\gg_i\mid X^*=-X\}=\kg\cap\gg_i.$$  
Also 
$$\ag_i=\{X=X^*\in\gg_i\mid [X,\Ec_i(X_0)]=0\}
\subseteq\{X=X^*\in\gg_i\mid [X,X_0]=0\}=\ag\cap\gg_i,$$ 
since $[\gg_i,\Ran(\1-\Ec_i)]=\{0\}.$ 
Finally, by Definitions \ref{D2}~and~\ref{iw_abstr1} we get 
$$\n_i=\gg_i\cap((\gg_i)_{\mathbb C})^{0,+}_{\ad(-\ie P(X_0))}
\subseteq\gg\cap(\gg_{\mathbb C})^{0,+}_{\ad(-\ie X_0)}\cap\gg_i
=\n\cap\gg_i,$$ 
where $(\bullet)_{\mathbb C}$ stands for the complexification of 
a Lie algebra.  
Consequently $K_i\subseteq K\cap G_i$, $A_i\subseteq A\cap G_i$, 
and $N_i\subseteq N\cap G_i$. 
(Here we use the fact that each of the groups $K_i$, $A_i$, and $N_i$ 
is connected since $G_i$ is connected 
and there exists a diffeomorphism $G_i\simeq K_i\times A_i\times N_i$.) 

Now we are going to use Lemma~\ref{top} 
with $\widetilde{S}_1=\widetilde{K}\times\widetilde{A}\times\widetilde{N}$, 
$\widetilde{S}_2=\widetilde{G}$, $S_1=K\times A\times N$, 
and $S_2=G$. 
The mapping $\widetilde{\m}$ is an open bijection since it 
is a diffeomorphism. 
On the other hand, 
$$\m(K\times A\times N)\supseteq
\bigcup_{i\in I}\m(K_i\times A_i\times N_i)
=\bigcup_{i\in I} G_i,$$
hence $\m(K\times A\times N)$ is a dense subset of $G$. 
(Note that $\bigcup\limits_{i\in I}G_i$ 
is dense in $G$ since $\bigcup\limits_{i\in I}\gg_i$ is dense in $\gg$ 
and $G$ is connected.) 
Thus Lemma~\ref{top} applies and shows that $\m(K\times A\times N)=G$, 
hence 
$\m\colon K\times A\times N\to G$ 
is a diffeomorphism. 

To complete the proof, use the inverse diffeomorphism 
$\m^{-1}\colon G\to K\times A\times N$. 
Since $\bigcup\limits_{i\in I}G_i$ 
is dense in $G$, it follows that 
$\bigcup\limits_{i\in I}K_i$ 
is dense in $K$, $\bigcup\limits_{i\in I}A_i$ 
is dense in $A$, and 
$\bigcup\limits_{i\in I}N_i$ 
is dense in $N$. 
Now the conclusion follows since $A_iN_i=N_iA_i$ 
and both groups $A_i$ and $N_i$ are simply connected for all $i\in I$.    
\end{proof}

\section{Classical Banach-Lie groups and 
their Lie algebras}\label{sectionClass}

In this section we introduce the Banach-Lie groups and 
Lie algebras whose Iwasawa decompositions 
will be investigated in 
Sections \ref{sectionA}, \ref{sectionB}, and \ref{sectionC} 
and we record a few auxiliary results that will be 
used in those sections.

\begin{definition}\label{complex_gr}
\normalfont
We denote by $\GL(\Hc)$ the group of all invertible bounded linear operators 
on the complex Hilbert space $\Hc$ and by 
 $\Ig$ an arbitrary norm ideal of $\Bc(\Hc)$. 
We define the following complex Banach-Lie groups and  
Banach-Lie algebras:  
\begin{itemize}
\item[{\rm(A)}] 
$\GL_{\Ig}(\Hc)=\GL(\Hc)\cap(\1+\Ig)$ 
with the Lie algebra 
$$\Lie(\GL_{\Ig}(\Hc)):=\gl_{\Ig}(\Hc):=\Ig;$$ 
\item[{\rm(B)}] 
$\OO_{\Ig}(\Hc):=\{g\in\GL_{\Ig}(\Hc)\mid g^{-1}=Jg^*J^{-1}\}$ 
with the Lie algebra 
$$\Lie(\OO_{\Ig}(\Hc))
:=\og_{\Ig}(\Hc):=\{x\in\Ig\mid x=-Jx^*J^{-1}\},$$ 
where $J\colon\Hc\to\Hc$ is a conjugation 
(i.e., $J$ a conjugate-linear isometry satisfying $J^2=\1$); 
\item[{\rm(C)}] 
$\Sp_{\Ig}(\Hc):=\{g\in\GL_{\Ig}(\Hc)\mid 
 g^{-1}=\widetilde{J}g^*\widetilde{J}^{-1}\}$ 
with the Lie algebra 
$$\Lie(\Sp_{\Ig}(\Hc))
:=\sp_{\Ig}(\Hc):=\{x\in\Ig\mid x=-\widetilde{J}x^*\widetilde{J}^{-1}\},$$ 
where $\widetilde{J}\colon\Hc\to\Hc$ is an anti-conjugation 
(i.e., $\widetilde{J}$ a conjugate-linear isometry satisfying 
$\widetilde{J}^2=-\1$).
\end{itemize}
We shall say that $\GL_{\Ig}(\Hc)$, $\OO_{\Ig}(\Hc)$, and $\Sp_{\Ig}(\Hc)$ 
are the \textit{classical complex Banach-Lie groups associated with  
the operator ideal $\Ig$}. 
Similarly, the corresponding Lie algebras 
are called the \textit{classical complex Banach-Lie algebras} 
(associated with $\Ig$). 

When no confusion can occur, we shall denote the groups 
$\GL_{\Ig}(\Hc)$, $\OO_{\Ig}(\Hc)$, and $\Sp_{\Ig}(\Hc)$ 
simply by 
$\GL_{\Ig}$, $\OO_{\Ig}$, and $\Sp_{\Ig}$, respectively, 
and we shall proceed similarly for the classical complex Lie algebras. 
\qed
\end{definition}

\begin{definition}\label{real_gr}
\normalfont
We shall use the notation of Definition~\ref{complex_gr} 
and define the following real Banach-Lie groups and Banach-Lie algebras:  
\begin{itemize}
\item[{\rm(AI)}] 
$\GL_{\Ig}(\Hc;{\mathbb R})=\{g\in\GL_{\Ig}(\Hc)\mid gJ=Jg\}$ 
with the Lie algebra 
$$\Lie(\GL_{\Ig}(\Hc;{\mathbb R})):=\gl_{\Ig}(\Hc;{\mathbb R})
:=\{x\in\Ig\mid xJ=Jx\},$$ 
where $J\colon\Hc\to\Hc$ is any conjugation on $\Hc$; 
\item[{\rm(AII)}]
$\GL_{\Ig}(\Hc;{\mathbb H})=\{g\in\GL_{\Ig}(\Hc)\mid 
 g\widetilde{J}=\widetilde{J}g\}$ 
with the Lie algebra 
$$\Lie(\GL_{\Ig}(\Hc;{\mathbb H})):=\gl_{\Ig}(\Hc;{\mathbb H})
:=\{x\in\Ig\mid x\widetilde{J}=\widetilde{J}x\},$$ 
where $\widetilde{J}\colon\Hc\to\Hc$ is any anti-conjugation on $\Hc$;
\item[{\rm(AIII)}] 
$\U_{\Ig}(\Hc_{+},\Hc_{-}):=\{g\in\GL_{\Ig}(\Hc)\mid g^*Vg=V\}$ 
with the Lie algebra 
$$\Lie(\U_{\Ig}(\Hc_{+},\Hc_{-})):=\ug_{\Ig}(\Hc_{+},\Hc_{-})
:=\{x\in\Jg\mid x^*V=-Vx\},$$
where $\Hc=\Hc_{+}\oplus\Hc_{-}$ and 
$V=\begin{pmatrix} \hfill 1 & \hfill 0 \\ \hfill 0 & \hfill -1\end{pmatrix}$ 
with respect to 
this orthogonal direct sum decomposition of $\Hc$; 
\item[{\rm(BI)}] 
$\OO_{\Ig}(\Hc_{+},\Hc_{-}):=\{g\in\GL_{\Ig}(\Hc)\mid g^{-1}=Jg^*J^{-1} 
\text{ and }g^*Vg=V\}$ 
with the Lie algebra 
$$\Lie(\OO_{\Ig}(\Hc_{+},\Hc_{-})):=\og_{\Ig}(\Hc_{+},\Hc_{-})
:=\{x\in\Jg\mid x=-Jx^*J^{-1}\text{ and }x^*V=-Vx\},$$
where $\Hc=\Hc_{+}\oplus\Hc_{-}$, 
$V=\begin{pmatrix} \hfill 1 & \hfill 0 \\ \hfill 0 & \hfill -1\end{pmatrix}$ 
with respect to 
this orthogonal direct sum decomposition of $\Hc$, 
and $J\colon\Hc\to\Hc$ is a conjugation on $\Hc$ 
such that $J(\Hc_{\pm})\subseteq\Hc_{\pm}$; 
\item[{\rm(BII)}]
$\OO^*_{\Ig}(\Hc):=\{g\in\GL_{\Ig}(\Hc)\mid g^{-1}=Jg^*J^{-1} 
\text{ and }g\widetilde{J}=\widetilde{J}g\}$ 
with the Lie algebra 
$$\Lie(\OO^*_{\Ig}(\Hc)):=\og^*_{\Ig}(\Hc)
:=\{x\in\Jg\mid x=-Jx^*J^{-1}\text{ and }x\widetilde{J}=\widetilde{J}x\},$$
where  $J\colon\Hc\to\Hc$ is a conjugation and 
$\widetilde{J}\colon\Hc\to\Hc$ is an anti-conjugation
such that $J\widetilde{J}=\widetilde{J}J$; 
\item[{\rm(CI)}] 
$\Sp_{\Ig}(\Hc;{\mathbb R}):=\{g\in\GL_{\Ig}(\Hc)\mid 
g^{-1}=\widetilde{J}g^*\widetilde{J}^{-1}\text{ and }gJ=Jg\}$ 
with the Lie algebra 
$$\sp_{\Ig}(\Hc;{\mathbb R}):=\{x\in\Ig\mid 
-x=\widetilde{J}x^*\widetilde{J}^{-1}\text{ and }xJ=Jx\},
$$
where $\widetilde{J}\colon\Hc\to\Hc$ is any anti-conjugation and  
$J\colon\Hc\to\Hc$ is any conjugation such that 
$J\widetilde{J}=\widetilde{J}J$; 
\item[{\rm(CII)}] 
$\Sp_{\Ig}(\Hc_{+},\Hc_{-}):=\{g\in\GL_{\Ig}(\Hc)\mid 
g^{-1}=\widetilde{J}g^*\widetilde{J}^{-1}\text{ and }g^*Vg=V\}$ 
with the Lie algebra 
$$\Lie(\Sp_{\Ig}(\Hc_{+},\Hc_{-})):=\sp_{\Ig}(\Hc_{+},\Hc_{-})
:=\{x\in\Jg\mid x=-\widetilde{J}x^*\widetilde{J}^{-1}\text{ and }x^*V=-Vx\},$$
where $\Hc=\Hc_{+}\oplus\Hc_{-}$, 
$V=\begin{pmatrix} \hfill 1 & \hfill 0 \\ \hfill 0 & \hfill -1\end{pmatrix}$ 
with respect to 
this orthogonal direct sum decomposition of $\Hc$, 
and $\widetilde{J}\colon\Hc\to\Hc$ is an anti-conjugation on $\Hc$ 
such that $\widetilde{J}(\Hc_{\pm})\subseteq\Hc_{\pm}$.
\end{itemize}
We say that  
$\GL_{\Ig}(\Hc;{\mathbb R})$, $\GL_{\Ig}(\Hc;{\mathbb H})$, 
$\U_{\Ig}(\Hc_{+},\Hc_{-})$, 
$\OO_{\Ig}(\Hc_{+},\Hc_{-})$, $\OO^*_{\Ig}(\Hc)$, 
$\Sp_{\Ig}(\Hc;{\mathbb R})$, and $\Sp_{\Ig}(\Hc_{+},\Hc_{-})$
are the \textit{classical real Banach-Lie groups associated with
the operator ideal $\Ig$}. 
Similarly, the corresponding Lie algebras 
are called the \textit{classical real Banach-Lie algebras} 
(associated with $\Ig$). 
\qed
\end{definition}

\begin{remark}\label{class_dlH}
\normalfont
The classical Banach-Lie groups and algebras of Definitions 
\ref{complex_gr}~and~\ref{real_gr} associated with 
the Schatten operator ideals $\Sg_p(\Hc)$ ($1\le p\le\infty$) 
were introduced in \cite{dlH72}, 
where it was conjectured that the connected $\1$-components of these groups 
have global Iwasawa decompositions in a natural sense
(see subsection~8.4 in Section~II.8 of \cite{dlH72}). 

We also note that as a by-product of 
the classification of the $L^*$-algebras 
(see for instance Theorems 7.18~and~7.19 in \cite{Bel06}), 
every (real or complex) topologically simple $L^*$-algebra 
is isomorphic to one of
the classical Banach-Lie algebras associated with 
the Hilbert-Schmidt ideal $\Jg=\Sg_2(\Hc)$.
\qed
\end{remark}

\begin{problem}\label{ideal}
\normalfont
In the setting of Definitions \ref{complex_gr}~and~\ref{real_gr}, 
the condition that $\Ig$ should be a \textit{norm} ideal is 
necessary in order to make the corresponding groups into 
smooth manifolds modeled on \textit{Banach} spaces 
(see for instance Proposition~9.28 in \cite{Bel06} or the beginning of 
the proof of Proposition~\ref{A} below).  
On the other hand, the classical ``Lie'' groups and Lie algebras 
can be defined with respect to any operator ideal, 
irrespective of whether it is endowed with a complete norm or not. 
And there exist lots of interesting operator ideals 
which do not support complete norms at all---see 
\cite{KW02}~and~\cite{KW06}.  

Thus it might prove important to study the Lie theoretic aspects 
of the classical groups and Lie algebras associated with arbitrary 
operator ideals, and perhaps to establish a bridge 
between the Lie theory and the commutator structure of 
operator ideals described in the papers \cite{DFWW04}~and~\cite{We05}. 
\qed
\end{problem}

We shall need the following generalization of 
Proposition~3 in \cite{Bal69}. 

\begin{lemma}\label{W0}
Let $\Hc$ be a complex separable infinite-dimensional Hilbert space,  
$\widetilde{J}\colon\Hc\to\Hc$ an anti-conjugation, 
and ${\mathbb K}\in\{{\mathbb R},{\mathbb C}\}$. 
Also let $Z\colon\Hc\to\Hc$ be a ${\mathbb K}$-linear continuous operator 
such that $Z\widetilde{J}=\widetilde{J}Z$ 
and $Z^2=z_0\1$ for some $z_0\in(0,\infty)$. 
Then there exists an orthonormal basis 
$\{\xi_l^{(\varepsilon)}
\}_{\stackrel{\scriptstyle l\in{\mathbb Z}\setminus\{0\}}
             {\scriptstyle\varepsilon\in\{\pm\sqrt{z_0}\}}}$ 
in the Hilbert space $\Hc$ over ${\mathbb K}$ 
such that 
$Z\xi_{\pm l}^{(\varepsilon)}=\varepsilon\xi_{\pm l}^{(\varepsilon)}$ 
and 
$\widetilde{J}\xi_{\pm l}^{(\varepsilon)}=\mp\xi_{\mp l}^{(\varepsilon)}$ 
whenever $\varepsilon\in\{\pm\sqrt{z_0}\}$
and 
$l=1,2,\dots$.  
\end{lemma}

\begin{proof}
Denote $\Hc^{(\varepsilon)}=\Ker(Z-\varepsilon)$ for 
$\varepsilon\in\{\pm\sqrt{z_0}\}$. 
Since $Z^2=z_0\1$, it follows that 
$\Hc=\Hc^{(\sqrt{z_0})}\oplus\Hc^{(-\sqrt{z_0})}$ 
as an orthogonal direct sum of ${\mathbb K}$-linear closed subspaces. 
On the other hand $\widetilde{J}Z=Z\widetilde{J}$, hence 
$\widetilde{J}\Hc^{(\varepsilon)}=\Hc^{(\varepsilon)}$ 
whenever $\varepsilon\in\{\pm\sqrt{z_0}\}$. 

Now let us keep $\varepsilon\in\{\pm\sqrt{z_0}\}$ fixed. 
We shall say that an orthonormal subset $\Sigma_0$ 
of the Hilbert space $\Hc^{(\varepsilon)}$ over ${\mathbb K}$ 
is a \textit{$\widetilde{J}$-set} if 
for each $x\in\Sigma_0$ we have 
$\{\widetilde{J}x,-\widetilde{J}x\}\cap\Sigma_0\ne\emptyset$. 
For every $x\in\Hc^{(\varepsilon)}$ with $\Vert x\Vert=1$ 
we have 
$(x\mid\widetilde{J}x)=-(\widetilde{J}^2x\mid\widetilde{J}x)
=-(x\mid\widetilde{J}x)$, 
whence $x\perp\widetilde{J}x$, 
so that $\{x,\widetilde{J}x\}$ is a $\widetilde{J}$-set. 
Then Zorn's lemma applies 
and shows that there exists 
a maximal $\widetilde{J}$-set $\Sigma^{(\varepsilon)}$ 
in the Hilbert space $\Hc^{(\varepsilon)}$ over ${\mathbb K}$. 

It is easy to see that $\Sigma^{(\varepsilon)}$ is 
actually an orthonormal basis 
in the Hilbert space $\Hc^{(\varepsilon)}$ over ${\mathbb K}$. 
In fact, let us assume that this is not the case and 
consider for instance the case ${\mathbb K}={\mathbb R}$. 
Then there exists $x_0\in\Hc^{(\varepsilon)}$ 
such that $\Vert x_0\Vert=1$ and $\Re(x\mid x_0)=0$ 
whenever $x\in\Sigma^{(\varepsilon)}$. 
Now for every $x\in\Sigma^{(\varepsilon)}$ we have 
$\Re(x\mid\widetilde{J}x_0)=-\Re(\widetilde{J}^2x\mid\widetilde{J}x_0)
=-\Re(\widetilde{J}x\mid x_0)=0$ 
since either $\widetilde{J}x\in\Sigma^{(\varepsilon)}$ or 
$-\widetilde{J}x\in\Sigma^{(\varepsilon)}$. 
It then follows that 
$\Sigma^{(\varepsilon)}\cup\{x_0,\widetilde{J}x_0\}$ is again a 
$\widetilde{J}$-set,  
thus contradicting the maximality of the 
$\widetilde{J}$-set~$\Sigma^{(\varepsilon)}$. 

Then let $\{\xi_l^{(\varepsilon)}\}_{l\ge1}$ be the set of all 
$x\in\Sigma^{(\varepsilon)}$ such that $-\widetilde{J}x\in\Sigma^{(\varepsilon)}$, 
and denote $\xi_{-l}^{(\varepsilon)}=-\widetilde{J}\xi_l^{(\varepsilon)}$ 
for $l=1,2,\dots$. 
Thus we get an orthonormal basis 
$\{\xi_l^{(\varepsilon)}
\}_{\stackrel{\scriptstyle l\in{\mathbb Z}\setminus\{0\}}
             {\scriptstyle\varepsilon\in\{\pm\sqrt{z_0}\}}}$ 
in the Hilbert space $\Hc$ over ${\mathbb K}$, 
satisfying the wished-for properties. 
\end{proof}

We shall also need the following version of 
Proposition~2 in \cite{Bal69}. 

\begin{lemma}\label{W-1}
Let $\Hc$ be a complex separable infinite-dimensional Hilbert space 
with a conjugation  
$J\colon\Hc\to\Hc$. 
Then the following assertions hold: 
\begin{itemize}
\item[{\rm(a)}] 
If $Z\in\Bc(\Hc)$ 
satisfies the conditions $ZJ=JZ$,   
$Z=Z^*=Z^{-1}$, and $\dim\Ker(Z-\1)=\dim\Ker(Z+\1)$, 
then there exists an orthonormal basis 
$\{\xi_l\}_{l\in{\mathbb Z}\setminus\{0\}}$ 
in the Hilbert space $\Hc$ 
such that 
$Z\xi_l=J\xi_l=\xi_{-l}$
whenever 
$l\in{\mathbb Z}\setminus\{0\}$. 
\item[{\rm(b)}] 
If $\widetilde{J}\colon\Hc\to\Hc$ is an anti-conjugation 
such that $J\widetilde{J}=\widetilde{J}J$, 
then there exists an orthonormal basis 
$\{\xi_l\}_{l\in{\mathbb Z}\setminus\{0\}}$ 
in the Hilbert space $\Hc$ 
such that $J\xi_l=\xi_{-l}$ whenever $l\in{\mathbb Z}\setminus\{0\}$,  
and $\widetilde{J}\xi_{\pm(2s-1)}=\xi_{\mp 2s}$ 
and $\widetilde{J}\xi_{\pm 2s}=-\xi_{\mp(2s-1)}$ 
for $s=1,2,\dots$. 
\end{itemize} 
\end{lemma}

\begin{proof}
We shall denote $\Hc_{\mathbb R}:=\{x\in\Hc\mid Jx=x\}$.  

(a) Let $\Hc_{\pm}=\Ker(Z\mp\1)$. 
Since $VJ=JV$, it follows that 
$\Hc_{\mathbb R}=(\Hc_{\mathbb R}\cap\Hc_{+})
\oplus(\Hc_{\mathbb R}\cap\Hc_{-})$, 
$\Hc_{\pm}=
(\Hc_{\mathbb R}\cap\Hc_{\pm})\oplus\ie(\Hc_{\mathbb R}\cap\Hc_{\pm})$,  
and $\dim_{\mathbb R}(\Hc_{\mathbb R}\cap\Hc_{\pm}) 
=\dim_{\mathbb C}(\Hc_{\pm})$ 
(see for instance Lemma~1 in \cite{Bal69}). 
Then there exist countable orthonormal bases in 
the real Hilbert spaces $\Hc_{\mathbb R}\cap\Hc_{\pm}$, 
which we denote by 
$\{x_{\pm l}\}_{l\ge1}$, 
respectively. 
Now define $\xi_{\pm l}=(x_l\pm\ie x_{-l})/\sqrt{2}$ 
for $l=1,2,\dots$. 
Then $\{\xi_l\}_{l\in{\mathbb Z}\setminus\{0\}}$ 
is the orthonormal basis in $\Hc$ which we were looking for. 

(b) Since $J\widetilde{J}=\widetilde{J}J$, 
we can use Lemma~\ref{W0} with $Z=J$ and ${\mathbb K}={\mathbb R}$ 
to get an orthonormal basis 
$\{x_l^{(\varepsilon)}
\}_{\stackrel{\scriptstyle l\in{\mathbb Z}\setminus\{0\}}
             {\scriptstyle\varepsilon\in\{\pm\}}}$ 
in the Hilbert space $\Hc$ over ${\mathbb R}$ 
such that 
$Jx_{\pm l}^{(\varepsilon)}=\varepsilon x_{\pm l}^{(\varepsilon)}$ 
and 
$\widetilde{J}x_{\pm l}^{(\varepsilon)}=\mp x_{\mp l}^{(\varepsilon)}$ 
whenever $\varepsilon\in\{\pm\}$
and 
$l=1,2,\dots$. 
 Let us denote $x_{l}=x_{l}^{+}$ for $l\in{\mathbb Z}\setminus\{0\}$. 
By using a bijection from $\{1,2,\dots\}$ 
onto the set of odd integer numbers, 
we can re-index the sequence $\{x_l\}_{l\in{\mathbb Z}\setminus\{0\}}$ 
such that $\widetilde{J}x_{\pm(2s-1)}=x_{\pm 2s}$ 
and $\widetilde{J}x_{\pm 2s}=-x_{\pm(2s-1)}$
for $s=1,2,\dots$. 
Now define just as above $\xi_{\pm l}=(x_l\pm\ie x_{-l})/\sqrt{2}$ 
for $l=1,2,\dots$, 
and thus we get an orthonormal basis in $\Hc$ 
satisfying the properties we wished for. 
\end{proof}

It will be convenient to have some special cases of 
Lemmas \ref{W0}~and~\ref{W-1} recorded as follows. 
(See Propositions 2~and~3 in \cite{Bal69} or \cite{Ne02b}.)

\begin{lemma}\label{W1}
Let $\Hc$ be a complex separable infinite-dimensional Hilbert space. 
\begin{itemize}
\item[{\rm(a)}] 
If $J\colon\Hc\to\Hc$ is a conjugation, 
then there exists an orthonormal basis 
$\{\xi_l\}_{l\in{\mathbb Z}\setminus\{0\}}$ 
in $\Hc$ such that 
$J\xi_l=\xi_{-l}$
whenever $l\in{\mathbb Z}\setminus\{0\}$. 
\item[{\rm(b)}] 
If $\widetilde{J}\colon\Hc\to\Hc$ is an anti-conjugation, 
then there exists an orthonormal basis 
$\{\widetilde{\xi}_l\}_{l\in{\mathbb Z}\setminus\{0\}}$ 
in $\Hc$ such that 
$\widetilde{J}\widetilde{\xi}_{\pm l}=\mp{\xi}_{\mp l}$ 
for $=1,2,\dots$. 
\end{itemize}
\end{lemma}

\begin{proof}
Assertion~(a) follows by our Lemma~\ref{W-1}(a) for a suitable choice of $Z$, 
while assertion~(b) follows by Lemma~\ref{W0} for $Z=\1$. 
\end{proof}

\begin{lemma}\label{W2}
Let $\Hc$ be a complex separable infinite-dimensional Hilbert space 
with an orthonormal basis $\{\xi_l\}_{l\in{\mathbb Z}\setminus\{0\}}$. 
Assume that $\{\alpha_l\}_{l\in{\mathbb Z}\setminus\{0\}}$ 
is a bounded family of real numbers  
and $\varepsilon\in{\mathbb R}$, and  
define 
the self-adjoint operator 
$$A=\sum\limits_{l\in{\mathbb Z}\setminus\{0\}}
 \alpha_l(\cdot\mid\xi_l)\xi_l\in\Bc(\Hc)$$ 
(where the sum is convergent in the strong operator topology). 
Then let $Z\colon\Hc\to\Hc$ be an ${\mathbb R}$-linear continuous operator 
such that $\Vert Zx\Vert=\Vert x\Vert$ whenever $x\in\Hc$ 
and satisfying either of the following conditions: 
\begin{itemize}
\item[{\rm(a)}] $Z^2=1$, $Z\xi_l=\xi_{-l}$ 
whenever $l\in{\mathbb Z}\setminus\{0\}$, 
and $Z$ is either ${\mathbb C}$-linear 
or ${\mathbb C}$-antilinear;  
\item[{\rm(b)}] 
$Z^2=-1$, $Z\xi_{\pm l}=\mp\xi_{\mp l}$ 
whenever $l=1,2,\dots$, 
and $Z$ is either ${\mathbb C}$-linear 
or ${\mathbb C}$-antilinear. 
\end{itemize}
Then we have $A=\varepsilon ZAZ^{-1}$ if and only if 
$\alpha_{-l}=\varepsilon\alpha_l$ 
whenever $l\in{\mathbb Z}\setminus\{0\}$. 
\end{lemma}

\begin{proof}
For every $l\in{\mathbb Z}\setminus\{0\}$ define 
the orthogonal projection $p_l=(\cdot\mid\xi_l)\xi_l$. 
Now assume that hypothesis~(a) is satisfied and $Z$ is antilinear. 
Then for every vector $\xi\in\Hc$ we have 
$$
Zp_lZ^{-1}\xi=Zp_lZ\xi=Z((Z\xi\mid\xi_l)\xi_l)
=(\xi_l\mid Z\xi)Z\xi_l=(Z^2\xi_l\mid Z\xi)Z\xi_l
=(\xi\mid Z\xi_l)Z\xi_l=(\xi\mid\xi_{-l})\xi_{-l}
=p_{-l}\xi,
$$
so that 
for all $l\in{\mathbb Z}\setminus\{0\}$ we have $Zp_lZ^{-1}=p_{-l}$. 
It is easy to see that the same conclusion holds if $Z$ were linear. 
Thence 
$$ZAZ^{-1}
=Z\Bigl(\sum\limits_{l\in{\mathbb Z}\setminus\{0\}}\alpha_lp_l\Bigr)Z^{-1}
=\sum\limits_{l\in{\mathbb Z}\setminus\{0\}}\alpha_lp_{-l},
$$
so that the equation $ZAZ^{-1}=\varepsilon A$ is equivalent to 
the fact that $\alpha_{-l}=\varepsilon\alpha_l$ 
whenever $l\in{\mathbb Z}\setminus\{0\}$. 

Now assume that hypothesis~(b) is satisfied and $Z$ is antilinear.
Then for all $l\in\{1,2,\dots\}$ and $\xi\in\Hc$ we have 
$$
Zp_{\pm l}Z^{-1}\xi
=-Zp_{\pm l}Z\xi
=-Z((Z\xi\mid\xi_{\pm l})\xi_{\pm l})
=-(\xi_{\pm l}\mid Z\xi)Z\xi_{\pm l}
=(Z^2\xi_{\pm l}\mid Z\xi)
  Z\xi_{\pm l}
=(\xi\mid Z\xi_{\pm l})Z\xi_{\pm l}.
$$
Consequently 
$Zp_lZ^{-1}=p_{-l}$ 
for every $l\in{\mathbb Z}\setminus\{0\}$, 
and we would obtain the same equality if $Z$ were linear. 
Then the wished-for conclusion follows just as above. 
\end{proof}

\begin{proposition}\label{A}
Let $\Hc$ be a complex Hilbert space and 
$\{e_i\}_{i\in I}$ a family of mutually orthogonal projections 
of rank~1 in $\Bc(\Hc)$ such that $\sum\limits_{i\in I}e_i=\1$. 
Now let $\Phi$ be a symmetric norming function, 
consider the corresponding Banach ideal $\Ig=\Sg_\Phi^{(0)}$, 
and define 
$$A:=\{g\in\GL_{\Ig}\mid g\ge 0;\; (\forall i\in I)\quad e_ig=ge_i\}
\text{ and }
\ag:=\{b\in\Ig\mid b=b^* ;\; (\forall i\in I)\quad e_ib=be_i\}.$$
Then $A$ is a simply connected Banach-Lie subgroup of $\GL_{\Ig}$ and 
$\Lie(A)=\ag$. 
\end{proposition}

\begin{proof} 
The proof is straightforward, using Proposition~4.4 in \cite{Bel06}. 
We omit the details. 
\end{proof}

\section{Iwasawa decompositions for groups of type A}~\label{sectionA}

Throughout this section we let 
$\Hc$ be a complex separable infinite-dimensional Hilbert space, 
$\Phi$ a mononormalizing symmetric norming function 
whose Boyd indices are non-trivial, 
and denote the corresponding separable norm ideal by 
$\Ig=\Sg_{\Phi}^{(0)}\subseteq\Bc(\Hc)$. 
We are going to apply the abstract methods developed in 
Section~\ref{sectionAbs} in order to construct 
global Iwasawa decompositions for classical groups of type~A 
associated with the operator ideal $\Ig$ 
(Definitions \ref{complex_gr}~and~\ref{real_gr}).

\addtocontents{toc}{\SkipTocEntry}
\subsection*{Complex groups of type A}

\begin{lemma}\label{matrices}
Let $\Hc_0$ be a complex Hilbert space and 
$S=-S^*\in\Bc(\Hc_0)$ 
with finite spectrum, say  $\sigma(S)=\{\lambda_1,\dots,\lambda_n\}$. 
For $k=1,\dots,n$ denote by $E_k\in\Bc(\Hc_0)$ 
the orthogonal projection onto $\Ker(S-\lambda_k)$. 
Then for every $Z\in\Bc(\Hc_0)$ 
we have 
$\int\limits_{\mathbb R}\ee^{t\cdot\ad\,S} Z\de\mu(t) 
=\sum\limits_{k=1}^nE_kZE_k$. 
\end{lemma}

\begin{proof}
We have $S=\sum\limits_{k=1}^n\lambda_k E_k$ 
and $E_kE_l=0$ whenever $k\ne l$, 
hence  
$\ee^{tS}=\sum\limits_{k=1}^n\ee^{t\lambda_k} E_k$ 
for arbitrary $t\in{\mathbb R}$. 
Thence  
\begin{equation}\label{block}
\ee^{t\cdot\ad\,S}Z=\ee^{tS}Z\ee^{-tS}
=\sum_{k,l=1}^n\ee^{t(\lambda_k-\lambda_l)} E_kZE_l.
\end{equation}
Now
$\int\limits_{\mathbb R}\ee^{t(\lambda_k-\lambda_l)}\de\mu(t)
=\int\limits_{\mathbb R}\ee^{(t+1)(\lambda_k-\lambda_l)}\de\mu(t)
=\ee^{\ie(\lambda_k-\lambda_l)}
 \int\limits_{\mathbb R}\ee^{t(\lambda_k-\lambda_l)}\de\mu(t) 
$
by the invariance property of $\mu$ (see Notation~\ref{D1}). 
If $k\ne l$ then $\lambda_k\ne \lambda_l$, whence 
$\int\limits_{\mathbb R}\ee^{t(\lambda_k-\lambda_l)}\de\mu(t)=0$. 
Then the wished-for equality follows by means of 
equation~\eqref{block}. 
\end{proof}

\begin{lemma}\label{approx}
Let $\Hc$ be a separable complex Hilbert space, 
$\Ig$ a separable norm ideal of $\Bc(\Hc)$, 
and $\{a_n\}_{n\ge0}$ a sequence of self-adjoint 
elements of $\Bc(\Hc)$ which is convergent 
to some $a\in\Bc(\Hc)$ in the strong operator topology. 
Then 
$\lim\limits_{n\to\infty}\Vert a_nx-ax\Vert_{\Jg}=
\lim\limits_{n\to\infty}\Vert xa_n-xa\Vert_{\Jg}=
\lim\limits_{n\to\infty}\Vert a_nxa_n-axa\Vert_{\Jg}=0$
whenever $x\in\Ig$. 
\end{lemma}

\begin{proof}
See Theorem~6.3 in Chapter~III of \cite{GK69}. 
\end{proof}

In the following statement, by \textit{triangular projection} 
associated with a self-adjoint operator we mean 
the triangular projection associated 
with its linearly ordered set of spectral projections; 
see \cite{GK70}, \cite{EL72}, \cite{Erd78}, and \cite{Pit88}.

\begin{proposition}\label{A_loc}
Let $X_0=X_0^*\in\gl_{\Ig}$ and denote 
$$\Lambda=\{\lambda\in{\mathbb R}\mid \dim\Ker(X_0-\lambda)>0\},$$
and for every $\lambda\in\Lambda$ 
let $E_\lambda\in\Bc(\Hc)$ 
be the orthogonal projection onto $\Ker(X_0-\lambda)$.  

Then $X_0$ is an Iwasawa quasi-regular element of $\gl_{\Ig}$
and  the following assertions hold: 
\begin{itemize}
\item[{\rm(1)}] 
The set $\Lambda$ is countable and $\{E_\lambda\}_{\lambda\in\Lambda}$ 
is a family of mutually orthogonal finite-rank projections 
satisfying $\sum\limits_{\lambda\in\Lambda}E_\lambda=\1$ 
in the strong operator topology. 
\item[{\rm(2)}] 
If $\Tc_{\Ig,X_0}\colon\gl_{\Ig}\to\gl_{\Ig}$ stands for the 
triangular projection associated with $X_0$, then 
\begin{eqnarray}
(\Tc_{\Ig,X_0})^2 &=& \Tc_{\Ig,X_0} 
  \nonumber\\
\Ran\Tc_{\Ig,X_0} &=& (\gl_{\Ig})_{\ad\,X_0}({\mathbb R}_{+}), 
   \nonumber\\
\Ran(\1-\Tc_{\Ig,X_0}) &\subseteq& (\gl_{\Ig})_{\ad\,X_0}(-{\mathbb R}_{+}), 
   \nonumber
\end{eqnarray}
\item[{\rm(3)}] 
Assume that we write the linear operators on $\Hc$ as infinite block matrices 
with respect to the partition of unity given by 
$\{E_\lambda\}_{\lambda\in\Lambda}$. 
Then for every $X\in\gl_{\Ig}$ the matrix of 
$\Dc_{\gl_{\Ig},\ad(-\ie X_0)}X$ can be constructed 
out of the matrix of $X$ by replacing all the off-diagonal blocks 
by zeros. 
\item[{\rm(4)}] 
If we denote 
$\ag_{\Ig,X_0}=\{X\in\gl_{\Ig}\mid X^*=X\text{ and }[X_0,X]=0\}$, 
then 
$$\ag_{\Ig,X_0}=\{X\in\gl_{\Ig}\mid X^*=X\text{ and }
(\forall\lambda\in\Lambda)\quad 
XE_\lambda\Hc\subseteq E_\lambda\Hc\}.
$$
\item[{\rm(5)}] 
We have 
$(\gl_{\Ig})_{\ad\,X_0}({\mathbb R}_{+})
=\Bigl\{X\in\gl_{\Ig}\mid(\forall\lambda\in\Lambda)\quad 
XE_\lambda\Hc\subseteq\spa
\Bigl(\bigcup\limits_{\lambda\le\beta\in\Lambda}E_\beta\Hc\Bigr)\Bigr\}$.
\item[{\rm(6)}] If we denote 
$\n_{\Ig,X_0}=(\gl_{\Ig})_{\ad(-\ie X_0)}^{0,+}$, 
then 
$$\n_{\Ig,X_0}=
\Bigl\{X\in\gl_{\Ig}\mid(\forall\lambda\in\Lambda)\quad 
XE_\lambda\Hc\subseteq\spa
\Bigl(\bigcup_{\lambda<\beta\in\Lambda}E_\beta\Hc\Bigr)\Bigr\}.$$
\item[{\rm(7)}] 
The Iwasawa decomposition of $\gl_{\Ig}$ associated with $X_0$ is 
$\gl_{\Ig}=\ug_{\Ig}\dotplus\ag_{\Ig,X_0}\dotplus\n_{\Ig,X_0}$.
\item[{\rm(8)}]
$X_0$ is an Iwasawa regular element of $\gl_{\Ig}$ if and only if 
it satisfies the condition 
\begin{equation}\label{mult1}
(\forall\lambda\in{\mathbb R})\quad \dim\Ker(X_0-\lambda)\le 1.
\end{equation}
\end{itemize}
\end{proposition}

\begin{proof} 
The fact that $X_0$ is an Iwasawa quasi-regular element will follow 
by Proposition~\ref{D4} as soon as we shall have proved assertion~(2). 

Assertion~(1) holds true since $X_0$ is a compact self-adjoint operator. 

We now prove assertion~(2). 
Since the Boyd indices of the symmetric norming function~$\Phi$ 
are nontrivial, it follows by Theorem~4.1 in \cite{Ara78} that 
the triangular projection defines a 
bounded linear idempotent operator $\Tc_{\Ig,X_0}\colon\gl_{\Ig}\to\gl_{\Ig}$. 
Let $\Fg$ stand for the ideal of finite-rank operators on $\Hc$. 
It is clear that 
$\Fg\cap\Ran\Tc_{\Ig,X_0}\subseteq(\gl_{\Ig})_{\ad\,X_0}({\mathbb R}_{+})$ 
and 
$\Fg\cap\Ran(\1-\Tc_{\Ig,X_0})\subseteq 
(\gl_{\Ig})_{\ad\,X_0}(-{\mathbb R}_{+})$. 
Since $\Fg$ is dense in $\Ig$ ($=\Sg_{\Phi}^{(0)}$) and 
$(\Tc_{\Ig,X_0})^2=\Tc_{\Ig,X_0} $, 
it easily follows that 
$\Fg\cap\Ran\Tc_{\Ig,X_0}$ is dense in $\Ran\Tc_{\Ig,X_0}$ 
and 
$\Fg\cap\Ran(\1-\Tc_{\Ig,X_0})$ is dense in $\Ran(\1-\Tc_{\Ig,X_0})$. 
Thence 
$$\Ran\Tc_{\Ig,X_0}\subseteq(\gl_{\Ig})_{\ad\,X_0}({\mathbb R}_{+}) 
\text{ and }
\Ran(\1-\Tc_{\Ig,X_0})\subseteq 
(\gl_{\Ig})_{\ad\,X_0}(-{\mathbb R}_{+}).$$
Thus, to complete the proof of assertion~(2), 
it remains to show that 
$\Ran\Tc_{\Ig,X_0}\supseteq(\gl_{\Ig})_{\ad\,X_0}({\mathbb R}_{+})$. 

To this end
let $X\in(\gl_{\Ig})_{\ad\,X_0}({\mathbb R}_{+})$ arbitrary. 
Then denote by $\{p_n\}_{n\ge0}$ any increasing sequence 
of finite-rank spectral projections of~$X_0$ such that 
$\lim\limits_{n\to\infty}p_n=\1$ in the strong operator topology. 
Then Lemma~\ref{approx} shows that 
$\lim\limits_{n\to\infty}\Vert p_nXp_n-X\Vert_{\Ig}=0$, 
hence it will be enough to check that $p_nXp_n\in\Ran\Tc_{\Ig,X_0}$ 
whenever $n\ge0$. 
For this purpose we recall that 
$$(\gl_{\Ig})_{\ad\,X_0}({\mathbb R}_{+})
=\Bigl\{Y\in\gl_{\Ig}\mid 
\limsup\limits_{t\to\infty} 
\frac{1}{t}\log\Vert(\exp(-t\cdot\ad\,X_0))Y\Vert\le0\Bigr\}
$$
(see for instance Remark~1.3 in \cite{Be01}). 
Since $X\in(\gl_{\Ig})_{\ad\,X_0}({\mathbb R}_{+})$ and 
$(\ad\,X_0)p_n=0$ whenever $n\ge0$, 
it then easily follows that 
$p_nXp_n\in(\gl_{\Ig})_{\ad\,X_0}({\mathbb R}_{+})$ 
for every $n\ge0$, and we have seen that this completes the 
proof of assertion~(2). 

To prove assertion~(3), let $\{p_n\}_{n\ge0}$ be 
a sequence of spectral projections of~$X_0$ as above. 
For every $n\ge0$ we have 
$(\ad(-\ie X_0))p_n=0$, 
whence it follows at once that 
$$(\forall n\ge0)\quad 
p_n(\Dc_{\gl_{\Ig},\ad(-\ie X_0)}X)p_n=\Dc_{\gl_{\Ig},\ad(-\ie X_0)}(p_nXp_n).
$$
Since $\Dc_{\gl_{\Ig},\ad(-\ie X_0)}\colon\gl_{\Ig}\to\gl_{\Ig}$ 
is a continuous operator and 
$\lim\limits_{n\to\infty}\Vert p_nYp_n-Y\Vert_{\Ig}=0$ 
for all $Y\in\Ig$, 
it then follows that it suffices to obtain the wished-for conclusion 
for the finite-rank operators $p_nXp_n$, 
where $n\ge0$ is arbitrary. 
And in this (finite-dimensional) case we just need to 
apply Lemma~\ref{matrices}. 

Assertion~(4) is a direct consequence of condition~\eqref{mult1}. 

Assertion~(5) follows at once by assertion~(2). 

Assertion~(6) follows by assertions (5)~and~(3). 

By using assertion~(2) along with Proposition~\ref{D4} 
we see that $X_0$ is an Iwasawa quasi-regular 
element of $\gl_{\Ig}$ and assertion~(7) holds.

Finally, we prove assertion~(8): 
$X_0$ is an Iwasawa regular element  
if and only if $[\ag_{\Ig,X_0},\ag_{\Ig,X_0}]=\{0\}$.  
And, by using assertion~(4), we see that the latter condition 
is equivalent to the fact that 
each eigenvalue of $X_0$ has the spectral multiplicity equal to~1, 
which is precisely condition~\eqref{mult1}. 
\end{proof}

We shall need the following extension of 
Lemma~5.2 in Chapter~VI of \cite{Hel01} 
to an infinite-dimensional setting. 

\begin{lemma}\label{L1}
Let $U$ be a real Banach-Lie group with the Lie algebra $\Lie(U)=\ug$, 
and assume that $\mg$ and $\hg$ are two closed subalgebras of $\ug$ such that 
the direct sum decomposition 
$\ug=\mg\dotplus\hg$
holds. 
Now let $M=\langle\exp_U(\mg)\rangle$ and 
$H=\langle\exp_U(\hg)\rangle$ be the corresponding subgroups of $U$ 
endowed with their natural structures of connected Banach-Lie groups 
such that $\Lie(M)=\mg$ and $\Lie(H)=\hg$. 
Then the multiplication map 
$\alpha\colon M\times H\to U$, $(m,h)\mapsto mh$, 
is smooth and has the property that for every $(m,h)\in M\times H$ 
the corresponding tangent map 
$T_{(m,h)}\colon T_{(m,h)}(M\times H)\to T_{mh}U$
is an isomorphism of Banach spaces. 
\end{lemma}

\begin{proof} 
The statement can be proved just as in the finite-dimensional case, 
so that we omit the details.
\end{proof}

\begin{theorem}\label{global}
Let $X_0=X_0^*\in\gl_{\Ig}$ satisfying 
condition~{\rm\eqref{mult1}}. 
Let 
$\Lambda=\{\lambda\in{\mathbb R}\mid \dim\Ker(X_0-\lambda)=1\}$,  
which is a linearly ordered set with respect to 
the reverse ordering of the real numbers, 
and for every $\lambda\in\Lambda$ pick 
$\xi_\lambda\in\Ker(X_0-\lambda)$ with $\Vert\xi_\lambda\Vert=1$.  
Now consider the Banach-Lie group 
$$G:=\GL_{\Ig}(\Hc)$$
and its subgroups 
$$\begin{aligned}
K:=&\{k\in G\mid k^*=k^{-1}\}, \\
A:=&\{a\in G\mid(\forall\lambda\in\Lambda)\quad 
a\xi_\lambda\in{\mathbb R}_{+}^*\xi_\lambda\}, 
 \text{ and }\\
N:=&\{n\in G\mid(\forall\lambda\in\Lambda)\quad 
n\xi_\lambda\in\xi_\lambda+\spa\{\xi_\beta\mid\beta<\lambda\}\}.
\end{aligned}$$
Then $K$, $A$, and $N$ are Banach-Lie subgroups of $G$, 
and the multiplication map 
$$\m\colon K\times A\times N\to G,\quad (k,a,n)\mapsto kan,$$
is a diffeomorphism. 
In addition, both subgroups $A$ and $N$ are simply connected and $AN=NA$. 
\end{theorem}

\begin{proof}
For every $\lambda\in\Lambda$ denote by 
$e_\lambda=(\cdot\mid \xi_\lambda)\xi_\lambda$ 
the orthogonal projection onto the one-dimensional subspace 
spanned by $\xi_\lambda$. 
Then it is easy to see that 
$$A=\{g\in\GL_{\Ig}\mid g\ge 0;\; (\forall\lambda\in\Lambda)\quad 
e_\lambda g=ge_\lambda\}, $$
hence Proposition~\ref{A}  
shows that $A$ is a Banach-Lie subgroup of $G$. 
Furthermore, the fact that $K$ is a Banach-Lie subgroup of $G$ 
follows e.g., by Proposition~9.28(ii) in \cite{Bel06}. 
As regards $N$, let us consider the Banach algebra 
$\Bg={\mathbb C}\1+\Ig$,
and note that 
$$N=\{n\in\Bg^\times\mid(\forall\lambda,\beta\in I,\lambda>\beta) 
\quad (n\xi_\lambda\mid\xi_\beta)=0;\; 
(\forall\lambda\in\Lambda)\quad (n\xi_\lambda\mid\xi_\lambda)=1;\; 
(\forall\lambda\in\Lambda)\quad np_\lambda=p_\lambda np_\lambda\},$$
hence $N$ is an algebraic subgroup of $\Bg^\times$, 
and thus it is a Banach-Lie subgroup of $\Bg^\times$ by 
the Harris-Kaup theorem (see \cite{HK77} or Theorem~4.13 in \cite{Bel06}). 
Since $N\subseteq G$ and $G$ is a Banach-Lie subgroup of $\Bg^\times$, 
it follows that $N$ is a Banach-Lie subgroup of $G$ as well. 

Now note that with the notation of Proposition~\ref{A_loc} 
we have 
$\Lie(G)=\gl_{\Ig}$, $\Lie(K)=\ug_{\Ig}$, $\Lie(A)=\ag_{\Ig,X_0}$, 
and $\Lie(N)=\n_{\Ig,X_0}$. 
By using Proposition~\ref{A}, 
it is easy to show that $B:=AN$ is a 
Banach-Lie subgroup of $G$ such that the multiplication mapping 
sets up a diffeomorphism $A\times N\to AN=B$. 
In addition, the Lie algebra of $B$ 
decomposes as $\Lie(B)=\ag_{\Ig,X_0}\dotplus\n_{\Ig,X_0}$. 
It then follows by Lemma~\ref{L1} and Proposition~\ref{A_loc}(7) 
that the multiplication mapping  
$\m\colon K\times A\times N\to G$, $(k,a,n)\mapsto kan$,  
is regular, in the sense that its tangent map 
$T_{(k,a,n)}\m\colon T_{(k,a,n)}(K\times A\times N)\to T_{kan}G$ 
is an isomorphism of Banach spaces for every 
$k\in K$, $a\in A$, and $n\in N$. 

Now define 
$$\psi\colon[0,1]\times (A\times N)\to A\times N,\quad 
\psi(t,a,n)=((1-t)a+t\1,\1+(1-t)(n-\1)).$$
Clearly $\1+(1-t)(n-\1)\in N$ whenever $n\in N$ and $t\in[0,1]$. 
On the other hand, for arbitrary $t\in(0,1]$ and $a\in A$ 
we have $(1-t)a+t\1\ge t\1$, 
hence $(1-t)a+t\1$ is a positive invertible operator; 
besides, 
$(1-t)a+t\1=a+t(\1-a)\in\1+\Sg_\Phi^{(0)}(\Hc)$, 
so that $(1-t)a+t\1\in A$. 
Consequently, the mapping $\psi$ is well defined. 
In addition, $\psi$ is continuous, 
$\psi(0,\cdot)=\id_{A\times N}$, 
and $\psi(1,\cdot)$ is a constant mapping of $A\times N$ into itself. 
Thus we easily see that both $A$ and $N$ are contractible topological spaces, 
and in particular they are simply connected. 

Since $K\cap AN=A\cap N=\{\1\}$, 
it is straightforward to prove that the mapping 
$\m\colon K\times A\times N\to G$ 
is injective. 
Since we have seen above that the mapping $\m$ is regular at every point, 
it then follows that $\m$ is a diffeomorphism of $K\times A\times N$ 
onto some open subset of $G$. 

To prove that the multiplication mapping $\m$ is actually surjective, 
let $g\in G$ arbitrary 
and consider the nest 
$\Pg=\{p_\lambda\}_{\lambda\in\Lambda}$, 
where $p_\lambda$ is the orthogonal projection onto the subspace  
$\Hc_\lambda=\spa\{\xi_\beta\mid \beta\le\lambda\}$ of $\Hc$ 
whenever $\lambda\in\Lambda$.  
Then  Corollary~\ref{A3}  shows that 
there exist a unitary element $w\in\1+\Ig$ and 
an element $b\in(\Alg\Pg)^\times\cap(\1+\Ig)$ 
such that $g=wb$. 
Now for every $\lambda\in\Lambda$ define 
$b_\lambda:=(b\xi_\lambda\mid\xi_\lambda)$. 
We have 
$\sup\limits_{\lambda\in\Lambda}\vert b_\lambda\vert\le\Vert v\Vert<\infty$, 
so that there exists an operator $d\in\Bc(\Hc)$ such that 
$$(d\xi_\lambda\mid\xi_\beta)
=\begin{cases}
b_\lambda & \text{ if }\lambda=\beta,\\
0 &\text{ if }\lambda\ne\beta. 
\end{cases}  
$$
Since $b\in(\Alg\Pg)^\times\cap(\1+\Ig)$, 
it follows at once that 
$d\in(\Alg\Pg)^\times\cap(\1+\Ig)$ as well, 
and then $d^{-1}b\in N$. 
Now let $d=u\vert d\vert$ be the polar decomposition of $d$. 
Then $u,\vert d\vert\in G$ by Lemma~5.1 in \cite{BR05} again, 
so that 
$g=wb=wu\vert d\vert(d^{-1}b)\in KAN$, 
and the proof ends. 
\end{proof} 

\addtocontents{toc}{\SkipTocEntry}
\subsection*{Real groups of type AI}

\begin{theorem}\label{AI_loc}
Let $J\colon\Hc\to\Hc$ be a conjugation 
and $\{\xi_l\}_{l\in{\mathbb Z}\setminus\{0\}}$ an orthonormal basis in $\Hc$ 
such that $J\xi_l=\xi_l$ whenever $l\in{\mathbb Z}\setminus\{0\}$. 
Pick a family of mutually different real numbers 
$\{\alpha_l\}_{l\in{\mathbb Z}\setminus\{0\}}$ 
such that  
$\lim\limits_{l\to\pm\infty}\alpha_l=0$
and 
$\Phi(\alpha_1,\alpha_{-1},\alpha_2,\alpha_{-2},\dots)<\infty$,  
and define the self-adjoint operator 
$$X_0=\sum\limits_{l\in{\mathbb Z}\setminus\{0\}}
 \alpha_l(\cdot\mid\xi_l)\xi_l\in\Bc(\Hc). $$ 
Then 
$X_0$ is an Iwasawa regular element of 
$\gl_{\Ig}(\Hc;{\mathbb R})$ 
and the Iwasawa decomposition of $\gl_{\Ig}(\Hc;{\mathbb R})$ 
associated with $X_0$ is 
\begin{equation}\label{AI_dec}
\gl_{\Ig}(\Hc;{\mathbb R})=
(\ug_{\Ig}\cap\gl_{\Ig}(\Hc;{\mathbb R}))
 \dotplus(\ag_{\Ig,X_0}\cap\gl_{\Ig}(\Hc;{\mathbb R}))
 \dotplus(\n_{\Ig,X_0}\cap\gl_{\Ig}(\Hc;{\mathbb R})) 
\end{equation}
(where $\ug_{\Ig}$, $\ag_{\Ig,X_0}$, and $\n_{\Ig,X_0}$ 
are the ones defined in {\rm Proposition~\ref{A_loc}}). 

Moreover, if $G$ stands for the connected $\1$-component of 
$\GL_{\Ig}(\Hc;{\mathbb R})$, 
then there exists a global Iwasawa decomposition 
$\m\colon K\times A\times N\to G$ corresponding to~\eqref{AI_dec}. 
In addition we have $AN=NA$, 
and both groups $A$ and $N$ are simply connected. 
\end{theorem}

\begin{proof}
The role of the orthonormal basis 
$\{\xi_l\}_{l\in{\mathbb Z}\setminus\{0\}}$ 
as in the statement can be played by any orthonormal basis 
in the real Hilbert space $\Hc_{\mathbb R}=\{x\in\Hc\mid Jx=x\}$. 
The conditions satisfied by the family 
$\{\alpha_l\}_{l\in{\mathbb Z}\setminus\{0\}}$ 
ensure that $X_0=X_0^*\in\Ig$. 
In addition, it is straightforward to check that actually 
$X_0\in\gl_{\Ig}(\Hc,{\mathbb R})$ 
(see for instance the proof of Lemma~\ref{W2}(a)). 

On the other hand, 
it follows by Proposition~\ref{A_loc} that $X_0$ is 
an Iwasawa quasi-regular element 
of $\gl_{\Ig}$ and the Iwasawa decomposition of $\gl_{\Ig}$ associated with 
$X_0$ is 
$$\gl_{\Ig}=\ug_{\Ig}\dotplus\ag_{\Ig,X_0}\dotplus\n_{\Ig,X_0}. $$
To obtain the Iwasawa decomposition 
asserted for $\gl_{\Ig}(\Hc,{\mathbb R})$ 
we are going to use
Corollary~\ref{D6} with 
$\widetilde{\gg}=\gl_{\Ig}$, $\gg_0=\gl_{\Ig}(\Hc;{\mathbb R})$, 
and $\widetilde{\Tc}=\Tc_{\Ig,X_0}\colon\gl_{\Ig}\to\gl_{\Ig}$ 
the triangular projection associated with $X_0$. 
To this end we have to prove that 
$\Tc_{\Ig,X_0}(\gl_{\Ig}(\Hc;{\mathbb R}))\subseteq 
\gl_{\Ig}(\Hc;{\mathbb R})$. 

In order to do so, we denote 
$\Hc_r=\spa\{\xi_1,\xi_{-1},\xi_2,\xi_{-2},\dots,\xi_r,\xi_{-r}\}$
and let $P_r\colon\Hc\to\Hc$ be the orthogonal projection onto $\Hc_r$ 
for $r=1,2,\dots$. 
Then for arbitrary $X\in\gl_{\Ig}(\Hc;{\mathbb R})\subseteq\Ig$ we have 
$\lim\limits_{r\to\infty}\Vert P_rXP_r-X\Vert_{\Ig}=0$  
by Lemma~\ref{approx}, hence 
$\lim\limits_{r\to\infty}
\Vert\Tc_{\Ig,X_0}(P_rXP_r)-\Tc_{\Ig,X_0}(X)\Vert_{\Ig}=0$. 
Thus it will be enough to show that 
$\Tc_{\Ig,X_0}(P_rXP_r)\in\gl_{\Ig}(\Hc;{\mathbb R})$ 
whenever $r\ge1$. 
And this follows by the restricted root space decomposition 
of the finite-dimensional real reductive Lie algebras 
$\gl(\Hc_r;{\mathbb R})\simeq\gl(r,{\mathbb R})$ for $r=1,2,\dots$. 

Thus $X_0$ is an Iwasawa quasi-regular element of 
$\gl_{\Ig}(\Hc;{\mathbb R})$. 
Since $X_0$ satisfies condition~\eqref{mult1} 
in Proposition~\ref{A_loc}, 
it follows that it is actually Iwasawa regular 
in $\gl_{\Ig}$, 
hence also in $\gl_{\Ig}(\Hc;{\mathbb R})$. 

To obtain the global Iwasawa decomposition, 
let us denote by $G$ the connected $\1$-component of 
$\GL_{\Ig}(\Hc;{\mathbb R})$. 
We are going to apply Proposition~\ref{lim_iw} 
with $\widetilde{G}$, $\widetilde{K}$, $\widetilde{A}$, and $\widetilde{N}$ 
as in Theorem~\ref{global}. 
To this end let $C$ be the connected $\1$-component of 
$\widetilde{C}\cap\GL_{\Ig}(\Hc;{\mathbb R})$ for $C\in\{K,A,N\}$. 
Then $C$ will be a connected Lie subgroup of $\widetilde{G}=\GL_{\Ig}$, 
since $\widetilde{C}\cap\GL_{\Ig}(\Hc;{\mathbb R})$ is a Lie subgroup 
of $\widetilde{G}$. 
(The latter property follows by the Harris-Kaup theorem if 
$C=K$ or $C=N$, and from Proposition~\ref{A} if $C=A$.) 
It is clear that 
$\Lie(G)=\gl_{\Ig}(\Hc;{\mathbb R})$, 
$\Lie(K)=\ug_{\Ig}\cap\gl_{\Ig}(\Hc;{\mathbb R})$, 
$\Lie(A)=\ag_{\Ig,X_0}\cap\gl_{\Ig}(\Hc;{\mathbb R})$, 
and 
$\Lie(N)=\n_{\Ig,X_0}\cap\gl_{\Ig}(\Hc;{\mathbb R})$, 
hence $\Lie(G)=\Lie(K)\dotplus\Lie(A)\dotplus\Lie(N)$. 
Next define 
$\Ec_r\colon\gl_{\Ig}(\Hc;{\mathbb R})\to\gl_{\Ig}(\Hc;{\mathbb R})$, 
$X\mapsto P_rXP_r$ 
for $r=1,2,\dots$. 
Now it is easy to see that Proposition~\ref{lim_iw} 
can be applied, and this completes the proof. 
\end{proof}

\addtocontents{toc}{\SkipTocEntry}
\subsection*{Real groups of type AII}

\begin{theorem}\label{AII_loc}
Let  $\widetilde{J}\colon\Hc\to\Hc$ be an anti-conjugation 
and $\{\widetilde{\xi}_l\}_{l\in{\mathbb Z}\setminus\{0\}}$ an 
orthonormal basis in $\Hc$ 
such that 
$\widetilde{J}\widetilde{\xi}_{\pm l}=\mp\widetilde{\xi}_{\mp l}$ 
for $l=1,2,\dots$. 
Pick a family of real numbers 
$\{\alpha_l\}_{l\in{\mathbb Z}\setminus\{0\}}$ 
such that the numbers 
$\{\alpha_l\}_{l\ge1}$ are mutually different and
$$\alpha_{-l}=\alpha_l\text{ for all }l\in{\mathbb Z}\setminus\{0\},
\quad  
 \lim\limits_{l\to\infty}\alpha_l=0,
\quad\text{and}\quad 
 \Phi(\alpha_1,\alpha_{-1},\alpha_2,\alpha_{-2},\dots)<\infty,$$
and define the self-adjoint operator 
$$X_0=\sum\limits_{l\in{\mathbb Z}\setminus\{0\}}
 \alpha_l(\cdot\mid\widetilde{\xi}_l)\widetilde{\xi}_l\in\Bc(\Hc). $$
Then  $X_0$ is an Iwasawa regular element of 
$\gl_{\Ig}(\Hc,{\mathbb H})$ and the Iwasawa decomposition of 
$\gl_{\Ig}(\Hc,{\mathbb H})$ associated with 
$X_0$ is 
\begin{equation}\label{AII_dec}
\gl_{\Ig}(\Hc;{\mathbb H})=
(\ug_{\Ig}\cap\gl_{\Ig}(\Hc;{\mathbb H}))
\dotplus(\ag_{\Ig,X_0}\cap\gl_{\Ig}(\Hc;{\mathbb H}))
\dotplus(\n_{\Ig,X_0}\cap\gl_{\Ig}(\Hc;{\mathbb H}))
\end{equation}
(where $\ug_{\Ig}$, $\ag_{\Ig,X_0}$, and $\n_{\Ig,X_0}$ 
are the ones defined in {\rm Proposition~\ref{A_loc}}). 

Moreover,  there exists a global Iwasawa decomposition 
$\m\colon K\times A\times N\to\GL_{\Ig}(\Hc;{\mathbb H})$ 
corresponding to~\eqref{AII_dec}. 
In addition we have $AN=NA$, 
and both groups $A$ and $N$ are simply connected. 
\end{theorem}

\begin{proof}
The orthonormal basis 
$\{\widetilde{\xi}_l\}_{l\in{\mathbb Z}\setminus\{0\}}$ 
as in the statement exists according to Lemma~\ref{W1}(a). 
The hypothesis on $\{\alpha_l\}_{l\in{\mathbb Z}\setminus\{0\}}$ 
implies that $X_0=X_0^*\in\Ig$, 
and then $X_0\in\gl_{\Ig}(\Hc;{\mathbb H})$ 
by Lemma~\ref{W2}(b). 

To see that $X_0$ is an Iwasawa quasi-regular element of 
$\gl_{\Ig}(\Hc;{\mathbb H})$ and the corresponding Iwasawa decomposition 
looks as asserted, 
one can proceed just as in the proof of Theorem~\ref{AI_loc}, 
this time using the orthogonal projection 
$\widetilde{P}_r\colon\Hc\to\Hc$ onto the subspace 
$\widetilde{\Hc}_r
=\spa\{\widetilde{\xi}_1,\widetilde{\xi}_{-1},
\widetilde{\xi}_2,\widetilde{\xi}_{-2},\dots,
\widetilde{\xi}_r,\widetilde{\xi}_{-r}\}
$
for $r=1,2,\dots$. 
We omit the details. 

It remains to check that $X_0$ is actually an Iwasawa regular element 
of $\gl_{\Ig}(\Hc;{\mathbb H})$. 
To this end denote 
$\widetilde{\Vc}_l={\mathbb C}\widetilde{\xi}_l+
{\mathbb C}\widetilde{\xi}_{-l}$ for $l=1,2,\dots$, 
and let $X=X^*\in\gl_{\Ig}(\Hc;{\mathbb H})$ 
with $[X,X_0]=0$. 
Since the real numbers $\{\alpha_l\}_{l\ge1}$ are mutually different 
and $[X,X_0]=0$ it follows that 
$X\widetilde{\Vc}_l\subseteq\widetilde{\Vc}_l$ whenever $l\ge1$. 
Now, since $X=X^*$, it follows that for each $l\ge1$ 
there exists an eigenvector $v_l\in\widetilde{\Vc}_l\setminus\{0\}$ of 
$X|_{\widetilde{\Vc}_l}$. 
Let $\gamma_l\in{\mathbb R}$ be the corresponding eigenvalue, 
so that $Xv_l=\gamma_lv_l$. 
On the other hand the anti-conjugation $\widetilde{J}$ satisfies 
$\widetilde{J}\widetilde{\Vc}_l\subseteq\widetilde{\Vc}_l$, 
hence $\widetilde{\Vc}_l$ has the natural structure of 
a quaternionic vector space. 
Since $\dim_{\mathbb C}\widetilde{\Vc}_l=2$ 
it follows that $\dim_{\mathbb H}\widetilde{\Vc}_l=1$, 
hence $\widetilde{\Vc}_l={\mathbb H}v_l$. 
Now the operator $X\in\gl_{\Ig}(\Hc;{\mathbb H})$ is ${\mathbb C}$-linear 
and $X\widetilde{J}=\widetilde{J}X$, 
hence for every $q\in{\mathbb H}$ we have 
$X(qv_l)=qXv_l=q\gamma_lv_l=\gamma_l(qv_l)$, 
so that $X\xi=\gamma_l\xi$ whenever $\xi\in\widetilde{\Vc}_l$. 

Since $\Hc=\bigoplus\limits_{l\ge1}\widetilde{\Vc}_l$, 
it then follows that $[X_1,X_2]=0$ whenever 
$X_j=X_j^*\in\gl_{\Ig}(\Hc;{\mathbb H})$ 
and $[X_j,X_0]=0$ for $j=1,2$. 

To prove the assertion on the global Iwasawa decomposition 
one can use Proposition~\ref{lim_iw} in a fashion similar 
to the one of the proof of Theorem~\ref{AI_loc}. 
\end{proof}

\addtocontents{toc}{\SkipTocEntry}
\subsection*{Real groups of type AIII}

\begin{theorem}\label{AIII_loc}
Assume that we have an orthogonal direct sum decomposition 
$\Hc=\Hc_{+}\oplus\Hc_{-}$ with $\dim\Hc_{+}=\dim\Hc_{-}$, 
and let $\{e_l^{\pm}\}_{l\ge1}$ be an orthonormal basis in $\Hc_{\pm}$. 
Then define $f_l^{\pm}=(e_l^{+}\pm e_l^{-})/\sqrt{2}$ 
whenever $l\ge1$. 
Pick a family of real numbers $\{\lambda_l\}_{l\ge1}$ 
such that 
$$\lambda_{j}\ne\pm\lambda_l\text{ if }j\ne l, 
\quad  
 \lim\limits_{l\to\infty}\lambda_l=0,
\quad\text{and}\quad 
 \Phi(\lambda_1,-\lambda_{1},\lambda_2,-\lambda_{2},\dots)<\infty,$$
and define the self-adjoint operator 
$$X_0:=\sum_{l\ge1}\lambda_l\bigl((\cdot\mid f_l^{+})f_l^{+}
-(\cdot\mid f_l^{-})f_l^{-}\bigr). $$
Then $X_0$ is an Iwasawa regular element of $\ug_{\Ig}(\Hc_{+},\Hc_{-})$ 
and the corresponding Iwasawa decomposition is 
\begin{equation}\label{AIII_dec}
\ug_{\Ig}(\Hc_{+},\Hc_{-})=
(\ug_{\Ig}\cap\ug_{\Ig}(\Hc_{+},\Hc_{-}))
\dotplus(\ag_{\Ig,X_0}\cap\ug_{\Ig}(\Hc_{+},\Hc_{-}))
\dotplus(\n_{\Ig,X_0}\cap\ug_{\Ig}(\Hc_{+},\Hc_{-})) 
\end{equation}
(where $\ug_{\Ig}$, $\ag_{\Ig,X_0}$, and $\n_{\Ig,X_0}$ 
are the ones defined in {\rm Proposition~\ref{A_loc}}).

Moreover, if $G$ stands for the connected $\1$-component of 
$\U_{\Ig}(\Hc_{+},\Hc_{-})$, 
then there exists a global Iwasawa decomposition 
$\m\colon K\times A\times N\to G$ corresponding to~\eqref{AI_dec}. 
In addition we have $AN=NA$, 
and both groups $A$ and $N$ are simply connected. 
\end{theorem}

\begin{proof}
To begin with, note that $\bigcup\limits_{l\ge1}\{f_l^{+},f_l^{-}\}$ 
is an orthonormal basis in $\Hc$. 
It then follows by the hypothesis on $\{\lambda_l\}_{l\ge1}$ 
along with Proposition~\ref{A_loc} 
that $X_0$ is an Iwasawa regular element of $\gl_{\Ig}$. 
We now show that actually $X_0\in\ug_{\Ig}(\Hc_{+},\Hc_{-})$. 
For this purpose denote 
$V=\begin{pmatrix} \hfill 1 & \hfill 0 \\
                   \hfill 0 & \hfill -1
   \end{pmatrix}$ as in Definition~\ref{real_gr}. 
Then for all $j\ge1$ we have $Ve_j^{\pm}=\pm e_j^{\pm}$, 
whence $Vf_j^{\pm}=f_j^{\mp}$. 
Thus for every $\xi\in\Hc$ we have 
$(\xi\mid f_j^{\pm})f_j^{\pm}=(\xi\mid Vf_j^{\mp})f_j^{\pm}
=V((V\xi\mid f_j^{\mp})f_j^{\mp})$. 
It then follows that $X_0=-VX_0V$, 
whence $X_0V=-VX_0$, and then $X_0\in\ug_{\Ig}(\Hc_{+},\Hc_{-})$. 

Now the wished-for conclusion will follow by Corollary~\ref{D6} 
as soon as we will have proved that the triangular projection 
$\Tc_{\Ig,X_0}\colon\Ig\to\Ig$ leaves 
$\ug_{\Ig}(\Hc_{+},\Hc_{-})$ invariant. 
To this end, for $r=1,2,\dots$
denote 
$$\Hc_r=\spa\{e_1^{+},e_1^{-},e_2^{+},e_2^{-},\dots,e_r^{+},e_r^{-}\}$$
and let $P_r\colon\Hc\to\Hc$ be the orthogonal projection onto $\Hc_r$. 
For arbitrary $X\in\ug_{\Ig}(\Hc_{+},\Hc_{-})\subseteq\Ig$ 
we have 
$\lim\limits_{r\to\infty}\Vert P_rXP_r-X\Vert_{\Ig}=0$ by Lemma~\ref{approx}, 
hence 
$\lim\limits_{r\to\infty}\Vert \Tc_{\Ig,X_0}(P_rXP_r)-
\Tc_{\Ig,X_0}(X)\Vert_{\Ig}=0$. 
Thus it will be enough to show that 
$\Tc_{\Ig,X_0}(P_rXP_r)\in\ug_{\Ig}(\Hc_{+},\Hc_{-})$ whenever $r\ge1$. 
And this follows by the restricted-root space decomposition 
of the finite-dimensional real reductive Lie algebras 
$\ug(\Hc_r\cap\Hc_{+},\Hc_r\cap\Hc_{-})
\simeq\ug(r,r)$ for $r=1,2,\dots$. 

To prove the assertion on the global Iwasawa decomposition 
one can use Proposition~\ref{lim_iw} in a fashion similar 
to the one of the proof of Theorem~\ref{AI_loc}. 
\end{proof}

\section{Iwasawa decompositions for groups of type B}\label{sectionB} 

As in Section~\ref{sectionA} we let 
$\Hc$ be a complex separable infinite-dimensional Hilbert space, 
$\Phi$ a mononormalizing symmetric norming function 
whose Boyd indices are non-trivial, 
and denote the corresponding separable norm ideal by 
$\Ig=\Sg_{\Phi}^{(0)}\subseteq\Bc(\Hc)$. 
We shall use the methods of  
Section~\ref{sectionAbs} to get 
global Iwasawa decompositions for classical groups of type~B
associated with the operator ideal $\Ig$.  

\addtocontents{toc}{\SkipTocEntry}
\subsection*{Complex groups of type B} 

\begin{theorem}\label{B_loc}
Let $J\colon\Hc\to\Hc$ be a conjugation 
and $\{\xi_l\}_{l\in{\mathbb Z}\setminus\{0\}}$ an orthonormal basis in $\Hc$ 
such that $J\xi_l=\xi_{-l}$ whenever $l\in{\mathbb Z}\setminus\{0\}$. 
Pick a family of mutually different real numbers 
$\{\alpha_l\}_{l\in{\mathbb Z}\setminus\{0\}}$ satisfying the conditions  
$$\alpha_{-l}=-\alpha_l\text{ for all }l\in{\mathbb Z}\setminus\{0\},
\quad  
 \lim\limits_{l\to\infty}\alpha_l=0,
\quad\text{and}\quad 
 \Phi(\alpha_1,\alpha_{-1},\alpha_2,\alpha_{-2},\dots)<\infty,$$ 
and define the self-adjoint operator 
$$X_0=\sum\limits_{l\in{\mathbb Z}\setminus\{0\}}
 \alpha_l(\cdot\mid\xi_l)\xi_l\in\Bc(\Hc). $$
Then  $X_0$ is an Iwasawa regular element of $\og_{\Ig}$ 
and the Iwasawa decomposition of $\og_{\Ig}$ associated with $X_0$ is 
\begin{equation}\label{B_dec}
\og_{\Ig}=
(\ug_{\Ig}\cap\og_{\Ig})\dotplus(\ag_{\Ig,X_0}\cap\og_{\Ig})
\dotplus(\n_{\Ig,X_0}\cap\og_{\Ig}) 
\end{equation}
(where $\ug_{\Ig}$, $\ag_{\Ig,X_0}$, and $\n_{\Ig,X_0}$ 
are the ones defined in {\rm Proposition~\ref{A_loc}}).

Moreover, if $G$ stands for the connected $\1$-component of 
$\OO_{\Ig}$, 
then there exists a global Iwasawa decomposition 
$\m\colon K\times A\times N\to G$ corresponding to~\eqref{B_dec}. 
In addition we have $AN=NA$, 
and both groups $A$ and $N$ are simply connected. 
\end{theorem}

\begin{proof}
Recall that the existence of the orthonormal basis 
$\{\xi_l\}_{l\in{\mathbb Z}\setminus\{0\}}$ 
as in the statement follows by Lemma~\ref{W1}(a). 
The conditions satisfied by the family 
$\{\alpha_l\}_{l\in{\mathbb Z}\setminus\{0\}}$ 
ensure that $X_0=X_0^*\in\Ig$. 
In addition, it follows by Lemma~\ref{W2}(a) that actually 
$X_0\in\og_{\Ig}$. 

On the other hand, since the real numbers in the family 
$\{\alpha_l\}_{l\in{\mathbb Z}\setminus\{0\}}$ are mutually different, 
it follows by Proposition~\ref{A_loc} that $X_0$ is an Iwasawa regular element 
of $\gl_{\Ig}$ and the Iwasawa decomposition of $\gl_{\Ig}$ associated with 
$X_0$ is 
$\gl_{\Ig}=\ug_{\Ig}\dotplus\ag_{\Ig,X_0}\dotplus\n_{\Ig,X_0}$. 
To obtain the conclusion we are going to use
Corollary~\ref{D5} for $\widetilde{\gg}=\gl_{\Ig}$, $\gg=\og_{\Ig}$, 
and $\widetilde{\Tc}=\Tc_{\Ig,X_0}\colon\gl_{\Ig}\to\gl_{\Ig}$. 
To this end it remains to prove that 
$\Tc_{\Ig,X_0}(\og_{\Ig})\subseteq\og_{\Ig}$. 

Denote 
$\Hc_r=\spa\{\xi_1,\xi_{-1},\xi_2,\xi_{-2},\dots,\xi_r,\xi_{-r}\}$
for $r\ge1$.
Also let $P_r\colon\Hc\to\Hc$ be the orthogonal projection onto $\Hc_r$ 
for $r=1,2,\dots$. 
Then for arbitrary $X\in\og_{\Ig}\subseteq\Ig$ we have 
$\lim\limits_{r\to\infty}\Vert P_rXP_r-X\Vert_{\Ig}=0$  
by Lemma~\ref{approx}, hence 
$\lim\limits_{r\to\infty}
\Vert\Tc_{\Ig,X_0}(P_rXP_r)-\Tc_{\Ig,X_0}(X)\Vert_{\Ig}=0$. 
Thus it will be enough to show that 
$\Tc_{\Ig,X_0}(P_rXP_r)\in\og_{\Ig}$ 
whenever $r\ge1$. 
And this follows by the restricted-root space decomposition 
of the finite-dimensional complex reductive Lie algebras 
$\og(\Hc_r)\simeq\og(r,{\mathbb C})$ for $r=1,2,\dots$. 

To prove the assertion on the global Iwasawa decomposition 
one can use Proposition~\ref{lim_iw} in a fashion similar 
to the one of the proof of Theorem~\ref{AI_loc}. 
\end{proof}

\addtocontents{toc}{\SkipTocEntry}
\subsection*{Real groups of type BI} 

\begin{theorem}\label{BI_loc}
Assume that $\Hc=\Hc_{+}\oplus\Hc_{-}$ with $\dim\Hc_{+}=\dim\Hc_{-}$
and let $J\colon\Hc\to\Hc$ be a conjugation 
such that $J(\Hc_{\pm})\subseteq\Hc_{\pm}$. 
Also let 
$V=\begin{pmatrix} \hfill 1 & \hfill 0 \\ \hfill 0 & \hfill -1\end{pmatrix}$ 
with respect to 
this orthogonal direct sum decomposition of $\Hc$. 
Then 
let $\{\xi_l\}_{l\in{\mathbb Z}\setminus\{0\}}$ an orthonormal basis in $\Hc$ 
such that $J\xi_l=V\xi_l=\xi_{-l}$ whenever $l\in{\mathbb Z}\setminus\{0\}$. 
Pick a family of mutually different real numbers 
$\{\alpha_l\}_{l\in{\mathbb Z}\setminus\{0\}}$ satisfying the conditions  
$$\alpha_{-l}=-\alpha_l\text{ for all }l\in{\mathbb Z}\setminus\{0\},
\quad  
 \lim\limits_{l\to\infty}\alpha_l=0,
\quad\text{and}\quad 
 \Phi(\alpha_1,\alpha_{-1},\alpha_2,\alpha_{-2},\dots)<\infty,$$ 
and define the self-adjoint operator 
$$X_0=\sum\limits_{l\in{\mathbb Z}\setminus\{0\}}
 \alpha_l(\cdot\mid\xi_l)\xi_l\in\Bc(\Hc). $$
Then  $X_0$ is an Iwasawa regular element of $\og_{\Ig}(\Hc_{+},\Hc_{-})$ 
and the corresponding Iwasawa decomposition 
is 
\begin{equation}\label{BI_dec}
\og_{\Ig}(\Hc_{+},\Hc_{-})=
(\ug_{\Ig}\cap\og_{\Ig}(\Hc_{+},\Hc_{-}))
\dotplus(\ag_{\Ig,X_0}\cap\og_{\Ig}(\Hc_{+},\Hc_{-}))
\dotplus(\n_{\Ig,X_0}\cap\og_{\Ig}(\Hc_{+},\Hc_{-})) 
\end{equation}
(where $\ug_{\Ig}$, $\ag_{\Ig,X_0}$, and $\n_{\Ig,X_0}$ 
are the ones defined in {\rm Proposition~\ref{A_loc}}).

Moreover if $G$ stands for the connected $\1$-component of 
$\OO_{\Ig}(\Hc_{+},\Hc_{-})$, 
then there exists a global Iwasawa decomposition 
$\m\colon K\times A\times N\to G$ corresponding to~\eqref{BI_dec}. 
In addition we have $AN=NA$, 
and both groups $A$ and $N$ are simply connected. 
\end{theorem}

\begin{proof}
The existence of the orthonormal basis 
$\{\xi_l\}_{l\in{\mathbb Z}\setminus\{0\}}$ 
as in the statement follows by Lemma~\ref{W-1}. 
The conditions satisfied by the family 
$\{\alpha_l\}_{l\in{\mathbb Z}\setminus\{0\}}$ 
ensure that $X_0=X_0^*\in\Ig$. 
In addition, it follows by Lemma~\ref{W2}(a) that actually 
$X_0\in\og_{\Ig}(\Hc_{+},\Hc_{-})$. 

On the other hand, since the real numbers in the family 
$\{\alpha_l\}_{l\in{\mathbb Z}\setminus\{0\}}$ are mutually different, 
it follows by Proposition~\ref{A_loc} that $X_0$ is an Iwasawa regular element 
of $\gl_{\Ig}$ and the Iwasawa decomposition of $\gl_{\Ig}$ associated with 
$X_0$ is 
$\gl_{\Ig}=\ug_{\Ig}\dotplus\ag_{\Ig,X_0}\dotplus\n_{\Ig,X_0}$. 
Thus the conclusion will follow by applying 
Corollary~\ref{D6} for $\widetilde{\gg}=\gl_{\Ig}$, 
$\gg_0=\og_{\Ig}(\Hc_{+},\Hc_{-})$, 
and $\widetilde{\Tc}=\Tc_{\Ig,X_0}\colon\gl_{\Ig}\to\gl_{\Ig}$. 
To this end it only remains to prove that 
$\Tc_{\Ig,X_0}(\og_{\Ig}(\Hc_{+},\Hc_{-}))
\subseteq\og_{\Ig}(\Hc_{+},\Hc_{-})$. 
In order to do so, we denote 
$$\Hc_r=\spa\{\xi_1,\xi_{-1},\xi_2,\xi_{-2},\dots,\xi_r,\xi_{-r}\}
\text{ whenever }r\ge1. $$
Also let $P_r\colon\Hc\to\Hc$ be the orthogonal projection onto $\Hc_r$ 
for $r=1,2,\dots$. 
Then for arbitrary $X\in\og_{\Ig}\subseteq\Ig$ we have 
$\lim\limits_{r\to\infty}\Vert P_rXP_r-X\Vert_{\Ig}=0$  
by Lemma~\ref{approx}, hence 
$$\lim\limits_{r\to\infty}
\Vert\Tc_{\Ig,X_0}(P_rXP_r)-\Tc_{\Ig,X_0}(X)\Vert_{\Ig}=0.$$
Thus it will be enough to show that 
$\Tc_{\Ig,X_0}(P_rXP_r)\in\og_{\Ig}(\Hc_{+},\Hc_{-})$ 
whenever $r\ge1$. 
And this follows by the restricted-root space decomposition 
of the finite-dimensional real reductive Lie algebras 
$\og(\Hc_r\cap\Hc_{+},\Hc_r\cap\Hc_{-})
\simeq\og(r,r)$ for $r=1,2,\dots$. 

To prove the assertion on the global Iwasawa decomposition 
one can use Proposition~\ref{lim_iw} in a fashion similar 
to the one of the proof of Theorem~\ref{AI_loc}. 
\end{proof}

\addtocontents{toc}{\SkipTocEntry}
\subsection*{Real groups of type BII}

\begin{theorem}\label{BII_loc}
Let $J\colon\Hc\to\Hc$ a conjugation and 
$\widetilde{J}\colon\Hc\to\Hc$ an anti-conjugation 
such that $J\widetilde{J}=\widetilde{J}J$. 
Then let 
$\{\xi_l\}_{l\in{\mathbb Z}\setminus\{0\}}$ 
be an orthonormal basis 
in the Hilbert space $\Hc$ 
such that $J\xi_l=\xi_{-l}$ whenever $l\in{\mathbb Z}\setminus\{0\}$,  
and $\widetilde{J}\xi_{\pm(2s-1)}=\xi_{\mp 2s}$ 
and $\widetilde{J}\xi_{\pm 2s}=-\xi_{\mp(2s-1)}$ 
for $s=1,2,\dots$. 

Pick a family of real numbers 
$\{\alpha_l\}_{l\in{\mathbb Z}\setminus\{0\}}$ 
such that 
$\alpha_{-l}=-\alpha_l$ for all $l\in{\mathbb Z}\setminus\{0\}$,
$\alpha_{2s-1}\ne\pm\alpha_{2t-1}$ whenever $s\ne t$ and $s,t\ge1$,   
$\alpha_{2s-1}=-\alpha_{2s}$ for $s=1,2,\dots$,
$\lim\limits_{l\to\infty}\alpha_l=0$,
and 
 $\Phi(\alpha_1,\alpha_{-1},\alpha_2,\alpha_{-2},\dots)<\infty$,  
and define the self-adjoint operator 
$$X_0=\sum\limits_{l\in{\mathbb Z}\setminus\{0\}}
 \alpha_l(\cdot\mid\xi_l)\xi_l\in\Bc(\Hc). $$
Then  $X_0$ is an Iwasawa regular element of $\og^*_{\Ig}(\Hc)$ 
and the corresponding Iwasawa decomposition 
is 
\begin{equation}\label{BII_dec}
\og^*_{\Ig}(\Hc)=
(\ug_{\Ig}\cap\og^*_{\Ig}(\Hc))
\dotplus(\ag_{\Ig,X_0}\cap\og^*_{\Ig}(\Hc))
\dotplus(\n_{\Ig,X_0}\cap\og^*_{\Ig}(\Hc)) 
\end{equation}
(where $\ug_{\Ig}$, $\ag_{\Ig,X_0}$, and $\n_{\Ig,X_0}$ 
are the ones defined in {\rm Proposition~\ref{A_loc}}).

Moreover, if $G$ stands for the connected $\1$-component of 
$\OO^*_{\Ig}(\Hc)$, 
then there exists a global Iwasawa decomposition 
$\m\colon K\times A\times N\to G$ corresponding to~\eqref{BII_dec}. 
In addition we have $AN=NA$, 
and both groups $A$ and $N$ are simply connected. 
\end{theorem}

\begin{proof}
The existence of the orthonormal basis 
$\{\xi_l\}_{l\in{\mathbb Z}\setminus\{0\}}$ 
as in the statement follows by Lemma~\ref{W-1}. 
Moreover, by Lemma~\ref{W2} we have $X_0\in\og_{\Ig}$. 
On the other hand, since 
$X_0=\sum\limits_{l\ge1}
 \alpha_l\bigl((\cdot\mid\xi_l)\xi_l-(\cdot\mid\xi_{-l})\xi_{-l}\bigr)$,
it follows (see the proof of Lemma~\ref{W2}(b)) that 
$$\widetilde{J}X_0\widetilde{J}^{-1}=\sum\limits_{l\ge1}
 \alpha_l\bigl((\cdot\mid\widetilde{J}\xi_l)\widetilde{J}\xi_l
-(\cdot\mid\widetilde{J}\xi_{-l})\widetilde{J}\xi_{-l}\bigr).$$
Now, since $\widetilde{J}\xi_{\pm(2s-1)}=\xi_{\mp 2s}$,  
$\widetilde{J}\xi_{\pm 2s}=-\xi_{\mp(2s-1)}$, and  
$\alpha_{2s-1}=-\alpha_{2s}$ for $s=1,2,\dots$, 
we see that $\widetilde{J}X_0\widetilde{J}^{-1}=X_0$. 
Thus $X_0\in\gl_{\Ig}(\Hc;{\mathbb H})\cap\og_{\Ig}=\og^*_{\Ig}(\Hc)$. 

By using the projections onto the subspaces 
$\Hc_s=\spa\{\xi_1,\xi_{-1},\xi_2,\xi_{-2},\dots,\xi_{2s},\xi_{-2s}\}$ 
and Corollary~\ref{D6} as in the proof of 
Proposition~\ref{BI_loc}, it then follows that 
$X_0$ is an Iwasawa quasi-regular element of $\og^*_{\Ig}(\Hc)$ 
and the corresponding Iwasawa decomposition looks as asserted. 
It remains to prove that $X_0$ is actually an Iwasawa regular 
element of $\og^*_{\Ig}(\Hc)$. 
And this fact can be obtained as in the proof 
of Proposition~\ref{AII_loc}. 
Specifically, denote $\Vc_s'={\mathbb C}\xi_{2s-1}+{\mathbb C}\xi_{-2s}$
$\Vc_s''={\mathbb C}\xi_{-(2s-1)}+{\mathbb C}\xi_{2s}$,  
and $\Vc_s=\Vc_s'\oplus\Vc_s''$ 
for $s=1,2,\dots$. 
Let $X=X^*\in\og^*_{\Ig}$ 
such that $[X,X_0]=0$. 
Since $\alpha_{2s-1}\ne\pm\alpha_{2t-1}$ whenever $s\ne t$ and $s,t\ge1$, 
it follows that $X$ leaves each of the subspaces $\Vc_s'$ and $\Vc_s''$ 
invariant for $s=1,2,\dots$. 
Since $X=X^*$, it follows that $X$ has 
eigenvectors $v_s'\in\widetilde{\Vc}_s'\setminus\{0\}$ 
and $v_s''\in\widetilde{\Vc}_s''\setminus\{0\}$. 
Let $\gamma_s',\gamma_s''\in{\mathbb R}$ be the corresponding eigenvalues, 
so that $Xv_s'=\gamma_s' v_s'$ and $Xv_s''=\gamma_s'' v_s''$. 
On the other hand the anti-conjugation $\widetilde{J}$ satisfies 
$\widetilde{J}\widetilde{\Vc}_s'\subseteq\widetilde{\Vc}_s'$ 
and 
$\widetilde{J}\widetilde{\Vc}_s''\subseteq\widetilde{\Vc}_s''$
hence both $\Vc_s'$ and $\Vc_s''$ have natural structures of 
quaternionic vector space. 
By counting dimensions, we get 
$\Vc_s'={\mathbb H}v_s'$ and $\Vc_s''={\mathbb H}v_s''$. 
The operator $X\in\gl_{\Ig}(\Hc;{\mathbb H})$ is ${\mathbb C}$-linear 
and $X\widetilde{J}=\widetilde{J}X$, 
hence for every $q\in{\mathbb H}$ we have 
$X(qv_s')=qXv_s'=\gamma_s'(qv_s')$, 
so that $X\xi=\gamma_s'\xi$ whenever $\xi\in\Vc_s'$. 
Similarly $X\xi=\gamma_s''\xi$ whenever $\xi\in\Vc_s''$.
Since $\Hc=\bigoplus\limits_{s\ge1}(\Vc_s'\oplus\Vc_s'')$, 
it then follows that $[X_1,X_2]=0$ whenever 
$X_j=X_j^*\in\og^*_{\Ig}(\Hc)$ 
and $[X_j,X_0]=0$ for $j=1,2$. 

To prove the assertion on the global Iwasawa decomposition 
one can use Proposition~\ref{lim_iw} in a fashion similar 
to the one of the proof of Theorem~\ref{AI_loc}. 
\end{proof}

\section{Iwasawa decompositions for groups of type C}~\label{sectionC}

Just as in Sections \ref{sectionA}~and~\ref{sectionB} we let 
$\Hc$ be a complex separable infinite-dimensional Hilbert space, 
$\Phi$ a mononormalizing symmetric norming function 
whose Boyd indices are non-trivial, 
and denote the corresponding separable norm ideal by 
$\Ig=\Sg_{\Phi}^{(0)}\subseteq\Bc(\Hc)$. 
As above, we shall use the methods of  
Section~\ref{sectionAbs} to get 
global Iwasawa decompositions for classical groups of type~C
associated with the operator ideal $\Ig$.  

\addtocontents{toc}{\SkipTocEntry}
\subsection*{Complex groups of type C}

\begin{theorem}\label{C_loc}
Let $\widetilde{J}\colon\Hc\to\Hc$ be an anti-conjugation 
and $\{\widetilde{\xi}_l\}_{l\in{\mathbb Z}\setminus\{0\}}$ an 
orthonormal basis in $\Hc$ 
such that 
$\widetilde{J}\widetilde{\xi}_{\pm l}=\mp\widetilde{\xi}_{\mp l}$ 
for $l=1,2,\dots$. 
Now pick a family of mutually different real numbers 
$\{\alpha_l\}_{l\in{\mathbb Z}\setminus\{0\}}$ satisfying the conditions 
$$\alpha_{-l}=-\alpha_l\text{ for all }l\in{\mathbb Z}\setminus\{0\},
\quad  
 \lim\limits_{l\to\infty}\alpha_l=0,
\quad\text{and}\quad 
 \Phi(\alpha_1,\alpha_{-1},\alpha_2,\alpha_{-2},\dots)<\infty,$$
and define the self-adjoint operator 
$$X_0=\sum\limits_{l\in{\mathbb Z}\setminus\{0\}}
 \alpha_l(\cdot\mid\widetilde{\xi}_l)\widetilde{\xi}_l\in\Bc(\Hc). $$
Then  $X_0$ is an Iwasawa regular element of 
$\sp_{\Ig}$ and the Iwasawa decomposition of $\sp_{\Ig}$ associated with 
$X_0$ is 
\begin{equation}\label{C_dec}
\sp_{\Ig}=
(\ug_{\Ig}\cap\sp_{\Ig})\dotplus(\ag_{\Ig,X_0}\cap\sp_{\Ig})
\dotplus(\n_{\Ig,X_0}\cap\sp_{\Ig}) 
\end{equation}
(where $\ug_{\Ig}$, $\ag_{\Ig,X_0}$, and $\n_{\Ig,X_0}$ 
are the ones defined in {\rm Proposition~\ref{A_loc}}).

Moreover there exists a global Iwasawa decomposition 
$\m\colon K\times A\times N\to \Sp_{\Ig}$ corresponding to~\eqref{C_dec}. 
In addition we have $AN=NA$, 
and both groups $A$ and $N$ are simply connected. 
\end{theorem}

\begin{proof}
One can proceed just as in the proof of Proposition~\ref{B_loc}, 
now using the orthogonal projection 
$\widetilde{P}_r\colon\Hc\to\Hc$ onto the subspace 
$\widetilde{\Hc}_r
=\spa\{\widetilde{\xi}_1,\widetilde{\xi}_{-1},
\widetilde{\xi}_2,\widetilde{\xi}_{-2},\dots,
\widetilde{\xi}_r,\widetilde{\xi}_{-r}\}
$
for $r=1,2,\dots$. 
We omit the details. 

To prove the assertion on the global Iwasawa decomposition 
one can use Proposition~\ref{lim_iw} in a fashion similar 
to the one of the proof of Theorem~\ref{AI_loc}. 
\end{proof}

\addtocontents{toc}{\SkipTocEntry}
\subsection*{Real groups of type CI}

\begin{theorem}\label{CI_loc}
Let $\widetilde{J}\colon\Hc\to\Hc$ be an anti-conjugation 
and $J\colon\Hc\to\Hc$ a conjugation such that 
$J\widetilde{J}=\widetilde{J}J$. 
Assume that 
$\{\widetilde{\xi}_l\}_{l\in{\mathbb Z}\setminus\{0\}}$ an 
orthonormal basis in $\Hc$ 
such that 
$\widetilde{J}\widetilde{\xi}_{\pm l}=\mp\widetilde{\xi}_{\mp l}$ 
and $J\widetilde{\xi}_{\pm l}=\widetilde{\xi}_{\pm l}$ 
for $l=1,2,\dots$. 
Now pick a family of mutually different real numbers 
$\{\alpha_l\}_{l\in{\mathbb Z}\setminus\{0\}}$ satisfying the conditions 
$$\alpha_{-l}=-\alpha_l\text{ for all }l\in{\mathbb Z}\setminus\{0\},
\quad  
 \lim\limits_{\alpha\to\infty}\alpha_l=0,
\quad\text{and}\quad 
 \Phi(\alpha_1,\alpha_{-1},\alpha_2,\alpha_{-2},\dots)<\infty,$$
and define the self-adjoint operator 
$$X_0=\sum\limits_{l\in{\mathbb Z}\setminus\{0\}}
 \alpha_l(\cdot\mid\widetilde{\xi}_l)\widetilde{\xi}_l\in\Bc(\Hc). $$
Then  $X_0$ is an Iwasawa regular element of 
$\sp_{\Ig}(\Hc;{\mathbb R})$ and the corresponding Iwasawa decomposition 
is 
\begin{equation}\label{CII_dec}
\sp_{\Ig}(\Hc;{\mathbb R})=
(\ug_{\Ig}\cap\sp_{\Ig}(\Hc;{\mathbb R}))
\dotplus(\ag_{\Ig,X_0}\cap\sp_{\Ig}(\Hc;{\mathbb R}))
\dotplus(\n_{\Ig,X_0}\cap\sp_{\Ig}(\Hc;{\mathbb R})) 
\end{equation}
(where $\ug_{\Ig}$, $\ag_{\Ig,X_0}$, and $\n_{\Ig,X_0}$ 
are the ones defined in {\rm Proposition~\ref{A_loc}}).

Moreover, if $G$ stands for the connected $\1$-component of 
$\Sp_{\Ig}(\Hc;{\mathbb R})$, 
then there exists a global Iwasawa decomposition 
$\m\colon K\times A\times N\to G$ corresponding to~\eqref{CI_dec}. 
In addition we have $AN=NA$, 
and both groups $A$ and $N$ are simply connected. 
\end{theorem}

\begin{proof}
The existence of an orthonormal basis 
as in the statement follows at once by Lemma~\ref{W0}. 

Let us prove that 
$X_0\in\sp_{\Ig}(\Hc;{\mathbb R})
=\sp_{\Ig}\cap\gl_{\Ig}(\Hc;{\mathbb R})$. 
In fact $X_0\in\sp_{\Ig}$ by Proposition~\ref{C_loc}. 
On the other hand, for all $l\in{\mathbb Z}\setminus\{0\}$ 
and $\eta\in\Hc$ we have 
$J\bigl((\eta\mid\widetilde{\xi}_l)\widetilde{\xi}_l\bigr)
=(\widetilde{\xi}_l\mid\eta)J\widetilde{\xi}_l
=(J\eta\mid J\widetilde{\xi}_l)\widetilde{\xi}_l
=(J\eta\mid\widetilde{\xi}_l)\widetilde{\xi}_l$, 
whence $JX_0=X_0J$, 
and thus $J\in\gl_{\Ig}(\Hc;{\mathbb R})$ as well. 

Moreover, just as in the proof of Theorem~\ref{AII_loc}, 
it follows that $X_0$ is an Iwasawa quasi-regular element of 
$\sp_{\Ig}(\Hc;{\mathbb R})$ and the corresponding Iwasawa decomposition 
looks as asserted. 
Finally, since $X_0$ is Iwasawa regular in $\sp_{\Ig}$ 
by Theorem~\ref{C_loc}, 
it follows that it is Iwasawa regular in $\sp_{\Ig}(\Hc;{\mathbb R})$ 
as well. 

To prove the assertion on the global Iwasawa decomposition 
one can use Proposition~\ref{lim_iw} in a fashion similar 
to the one of the proof of Theorem~\ref{AI_loc}. 
\end{proof}

\addtocontents{toc}{\SkipTocEntry}
\subsection*{Real groups of type CII}

\begin{theorem}\label{CII_loc}
Assume that we have an orthogonal direct sum decomposition 
$\Hc=\Hc_{+}\oplus\Hc_{-}$ with $\dim\Hc_{+}=\dim\Hc_{-}$, 
and let $\widetilde{J}\colon\Hc\to\Hc$ be an anti-conjugation 
such that $\widetilde{J}(\Hc_{\pm})\subseteq\Hc_{\pm}$. 
Also let 
$V=\begin{pmatrix} \hfill 1 & \hfill 0 \\ \hfill 0 & \hfill -1\end{pmatrix}$ 
with respect to 
this orthogonal direct sum decomposition of $\Hc$. 
Now let $\bigcup\limits_{l\in{\mathbb Z}\setminus\{0\}}
\{e_l^{+},e_l^{-}\}$ 
be an orthonormal basis in $\Hc$ such that 
$\widetilde{J}e^{\varepsilon}_{\pm l}=\mp e^{\varepsilon}_{\mp l}$ 
and $Ve^{\varepsilon}_{\pm l}=\varepsilon e^{\varepsilon}_{\pm l}$ 
whenever $\varepsilon\in\{+,-\}$ and $l=1,2,\dots$. 
Then define $f_l^{\pm}=(e_l^{+}\pm e_l^{-})/\sqrt{2}$ 
for all $l\in{\mathbb Z}\setminus\{0\}$. 
Pick a family of mutually different real numbers 
$\{\lambda_l\}_{l\in{\mathbb Z}\setminus\{0\}}$ 
such that 
$$\lambda_{-l}=-\lambda_l\text{ whenever }l\in{\mathbb Z}\setminus\{0\}, 
\quad
\lim\limits_{l\to\infty}\lambda_l=0,
\quad\text{and}\quad 
 \Phi(\lambda_1,-\lambda_1,\lambda_2,-\lambda_2,\dots)<\infty,$$
and define the self-adjoint operator 
$$X_0:=\sum_{l\in{\mathbb Z}\setminus\{0\}}
\lambda_l\bigl((\cdot\mid f_l^{+})f_l^{+}
-(\cdot\mid f_l^{-})f_l^{-}\bigr). $$
Then $X_0$ is an Iwasawa regular element of $\sp_{\Ig}(\Hc_{+},\Hc_{-})$ 
and the corresponding Iwasawa decomposition is 
\begin{equation}\label{CI_dec}
\sp_{\Ig}(\Hc_{+},\Hc_{-})=
(\ug_{\Ig}\cap\sp_{\Ig}(\Hc_{+},\Hc_{-}))
\dotplus(\ag_{\Ig,X_0}\cap\sp_{\Ig}(\Hc_{+},\Hc_{-}))
\dotplus(\n_{\Ig,X_0}\cap\sp_{\Ig}(\Hc_{+},\Hc_{-})) 
\end{equation}
(where $\ug_{\Ig}$, $\ag_{\Ig,X_0}$, and $\n_{\Ig,X_0}$ 
are the ones defined in {\rm Proposition~\ref{A_loc}}).

Moreover, if $G$ stands for the connected $\1$-component of 
$\Sp_{\Ig}(\Hc_{+},\Hc_{-})$, 
then there exists a global Iwasawa decomposition 
$\m\colon K\times A\times N\to G$ corresponding to~\eqref{CII_dec}. 
In addition we have $AN=NA$, 
and both groups $A$ and $N$ are simply connected. 
\end{theorem}

\begin{proof}
The existence of the orthonormal basis 
$\bigcup\limits_{l\in{\mathbb Z}\setminus\{0\}}
\{e_l^{+},e_l^{-}\}$ follows by Lemma~\ref{W0}. 
Just as in the proof of Proposition~\ref{AIII_loc} 
we see that $X_0\in\ug_{\Ig}(\Hc_{+},\Hc_{-})$. 
On the other hand, Lemma~\ref{W2}(b) 
shows that $X_0\in\sp_{\Ig}$, 
and thus $X_0\in\sp_{\Ig}\cap\ug_{\Ig}(\Hc_{+},\Hc_{-})
=\sp_{\Ig}(\Hc_{+},\Hc_{-})$. 
Then, by using Corollary~\ref{D6} along with 
the orthogonal projections on the 
subspaces $\Hc_r=\spa\Bigl(\bigcup\limits_{-r\le l\le r}
\{e_l^{+},e_l^{-}\}\Bigr)$ for $r=1,2,\dots$, 
one can prove that $X_0$ is an Iwasawa quasi-regular element 
of $\sp_{\Ig}(\Hc_{+},\Hc_{-})$ 
and the corresponding Iwasawa decomposition 
looks as asserted. 
(See 
the proof of Proposition~\ref{AIII_loc} for some more details.) 

Now it remains to show that $X_0$ is actually an 
Iwasawa regular element 
of $\sp_{\Ig}(\Hc_{+},\Hc_{-})$. 
To this end denote 
$\Vc_j^0={\mathbb C}f_j^{+}+{\mathbb C}f_{-j}^{-}$,  
$\Vc_j^1={\mathbb C}f_j^{-}+{\mathbb C}f_{-j}^{+}$, 
and $\Vc_j=\Vc_j^0\oplus\Vc_j^1$  
for $j\ge 1$. 
Then let $X=X^*\in\sp_{\Ig}(\Hc_{+},\Hc_{-})$ 
such that $[X,X_0]=0$. 
We have  $\lambda_{-l}=-\lambda_l$ whenever 
$l\in{\mathbb Z}\setminus\{0\}$, hence
$$X_0:=\sum_{j\ge1}
\lambda_j\bigl((\cdot\mid f_j^{+})f_j^{+}
+(\cdot\mid f_{-j}^{-})f_{-j}^{-}\bigr)
-\lambda_j\bigl((\cdot\mid f_j^{-})f_j^{-}
+(\cdot\mid f_{-j}^{+})f_{-j}^{+}\bigr). $$
Since the real numbers $\{\lambda_j\}_{j\ge1}$ are mutually 
different and $[X,X_0]=0$, 
it follows that $X$ leaves both the subspaces $\Vc_j^0$ and $\Vc_j^1$ 
invariant whenever $j=1,2,\dots$. 
Now let us keep $j\in\{1,2,\dots\}$ and $\varepsilon\in\{0,1\}$ fixed. 
Since $X=X^*$, there exist $x_0\in\Vc_j^\varepsilon\setminus\{0\}$ and 
$t_0\in{\mathbb R}$ 
such that $Xx_0=t_0x_0$. 
On the other hand, since $V\widetilde{J}=\widetilde{J}V$, 
it follows directly that $\widetilde{J}_1:=V\widetilde{J}$ 
is an anti-conjugation on $\Hc$. 
In addition, since $Vf_{\pm j}^\varepsilon=f_{\pm j}^{-\varepsilon}$ 
and $\widetilde{J}f_{\pm j}^\varepsilon=\mp f_{\mp j}^\varepsilon$, 
it follows that the linear subspace $\Vc_j^\varepsilon$ 
is invariant under the anti-conjugation $\widetilde{J}_1$. 
Let us endow $\Vc_j^\varepsilon$ with 
the corresponding quaternionic structure. 
Since $\dim_{\mathbb C}\Vc_j^\varepsilon=2$, 
it follows that $\dim_{\mathbb H}\Vc_j^\varepsilon=1$, 
and thus $\Vc_j^\varepsilon={\mathbb H}x_0$. 
On the other hand $X=X^*\in\sp_{\Ig}(\Hc_{+},\Hc_{-})$, 
hence $XV=-VX$ and $X\widetilde{J}=-\widetilde{J}X$, 
whence $X\widetilde{J}_1=\widetilde{J}_1X$. 
Thus $X$ is an ${\mathbb H}$-linear operator 
with respect to the quaternionic structure defined by the anti-conjugation 
$\widetilde{J}_1$. 
Now, since $\Vc_j^\varepsilon={\mathbb H}x_0$ 
and $Xx_0=t_0x_0$, 
it follows that the restriction of $X$ to $\Vc_j^\varepsilon$ 
is given by the multiplication by the real number $t_0$. 
Since $\Hc=\bigoplus\limits_{j\ge1}(\Vc_j^0\oplus\Vc_j^1)$, 
it thus follows that the operators in 
$\ag_{\Ig,X_0}\cap\sp_{\Ig}(\Hc_{+},\Hc_{-})$ 
commute pairwise, and this completes the proof. 

To prove the assertion on the global Iwasawa decomposition 
one can use Proposition~\ref{lim_iw} in a fashion similar 
to the one of the proof of Theorem~\ref{AI_loc}. 
\end{proof}

\section{Group decompositions for covering groups}\label{sectionCov}

The aim of this short section is to 
show that the Iwasawa decompositions constructed in 
Sections \ref{sectionA}, \ref{sectionB}, and \ref{sectionC} 
can be lifted to any covering groups. 
We refer to \cite{dlH72} and \cite{Ne02b} 
for information on the homotopy groups of the classical Banach-Lie groups 
associated with the Schatten ideals. 
It is easy to see that the correspnding description 
of homotopy groups actually holds true 
for the classical Banach-Lie groups associated with 
any separable norm ideal. 

\begin{proposition}\label{lift_iw}
Let $G$ be a connected Banach-Lie group, and $K$, $A$, and $N$ 
connected Banach-Lie subgroups of $G$ such that the multiplication map
$\m\colon K\times A\times N\to G$
is a diffeomorphism. 
In addition, assume that $A$ and $N$ are simply connected 
and $AN=NA$. 

Now assume that we have a connected Banach-Lie group $\widetilde{G}$ 
with a covering homomorphism 
$e\colon\widetilde{G}\to G$, and define 
$\widetilde{K}:=e^{-1}(K)$, $\widetilde{A}:=e^{-1}(A)$, 
and $\widetilde{N}:=e^{-1}(N)$. 
Then $\widetilde{K}$, $\widetilde{A}$, and $\widetilde{N}$ 
are connected Banach-Lie subgroups of $\widetilde{G}$ and the multiplication 
map 
$\widetilde{\m}\colon\widetilde{K}\times\widetilde{A}\times\widetilde{N}\to 
\widetilde{G}$
is a diffeomorphism. 
\end{proposition}

\begin{proof} 
The proof can be achieved by using straightforward infinite-dimensional versions 
of some standard ideas from the theory of Iwawsawa decompositions 
of reductive groups 
(specifically, see the proofs of Theorem 6.31~and~6.46 in \cite{Kna96}). 
We omit the details. 
\end{proof}

\begin{corollary}\label{covA} 
Let $\Hc$ be a complex separable infinite-dimensional Hilbert space, 
$\Phi$ a mononormalizing symmetric norming function 
whose Boyd indices are non-trivial, 
and denote the corresponding separable norm ideal by 
$\Ig=\Sg_{\Phi}^{(0)}\subseteq\Bc(\Hc)$. 
Then let $\m\colon K\times A\times N\to G$ be the global Iwasawa decomposition 
given by any of 
{\rm Theorems \ref{global}, \ref{AI_loc}, \ref{AII_loc}, \ref{AIII_loc}, 
\ref{BI_loc}, \ref{BII_loc}, 
\ref{C_loc}, \ref{CI_loc}~\textit{and}~\ref{CII_loc}} 
for the connected $\1$-components of real or complex classical 
Banach-Lie groups.   
Now denote by $e\colon\widehat{G}\to G$ any covering group of $G$. 
If we define $\widehat{K}:=e^{-1}(K)$, $\widehat{A}:=e^{-1}(A)$, 
and $\widehat{N}:=e^{-1}(N)$,  
then $\widehat{K}$, $\widehat{A}$, and $\widehat{N}$ 
are connected Banach-Lie subgroups of $\widetilde{G}$ and the multiplication 
map 
$\widehat{\m}\colon\widehat{K}\times\widehat{A}\times\widehat{N}\to 
\widehat{G}$
is a diffeomorphism. 
\end{corollary}

\begin{proof}
Use Proposition~\ref{lift_iw}. 
\end{proof}

\appendix 

\section{Auxiliary facts on operator ideals}\label{App}

In this appendix we record some facts on operator ideals, 
stating them under versions appropriate for use in the main body 
of the present paper. 
We refer to \cite{GK69}, \cite{GK70}, \cite{Erd72}, \cite{EL72}, \cite{Erd78}, 
\cite{KW02}, \cite{We05}, \cite{DFWW04}, \cite{KW06},  and \cite{Bel06} for 
various special topics involving symmetric norm ideals 
related to the circle of ideas discussed here. 

Let $\Hc$ be a complex Hilbert space
and $\Pg$ a maximal nest in $\Bc(\Hc)$. 
That is, $\Pg$ is a maximal linearly ordered set 
of orthogonal projections on $\Hc$. 
Then we denote 
$\Alg\Pg:=\{b\in\Bc(\Hc)\mid (\forall p\in\Pg)\quad bp=pbp\}$
(the \textit{nest algebra} associated with $\Pg$). 

In the following statement we need the notion of Boyd indices 
as used in \cite{Ara78} (see also subsections 2.17--19 in \cite{DFWW04}). 

\begin{theorem}\label{A2}
Assume that $\Hc$ is a complex separable Hilbert space and 
$\Pg$ is a maximal nest in $\Bc(\Hc)$. 
Let $\Phi$ be a symmetric norming function whose Boyd indices are nontrivial 
and denote $\Ig=\Sg_\Phi^{(0)}$. 
Then for every $a\in\GL_{\Ig}(\Hc)$ such that $0\le a$ 
there exist uniquely determined operators $d\in\GL_{\Ig}(\Hc)$ 
and $r\in\Ig$ satisfying the following conditions: 
\begin{itemize}
\item[$\bullet$] $0\le d\in\GL_{\Ig}(\Hc)\cap\Alg\Pg$; 
\item[$\bullet$] $r\in\Ig\cap\Alg\Pg$ and 
    the spectrum of $r$ is equal to $\{0\}$; 
\item[$\bullet$] $a=(\1+r^*)d(\1+r)$. 
\end{itemize}
\end{theorem}

\begin{proof}
Theorem~4.1 in \cite{Ara78} and Lemma~4.3 in \cite{EL72} 
show that Theorem~4.2 in \cite{Erd72} 
(or Theorem~6.2 in Chapter~IV of \cite{GK70}) applies 
for the operator ideals $\Sg_{\text{I}}=\Sg_{\text{II}}=\Sg_\Phi^{(0)}$. 
\end{proof}

\begin{corollary}\label{A3}
Let $\Pg$, $\Phi$, and $\Ig$ be as in {\rm Theorem~\ref{A2}}. 
Then for every $g\in\GL_{\Ig}(\Hc)$ there exist the operators 
$b\in\GL_{\Ig}\cap\Alg\Pg$ and $u\in\U_{\Ig}(\Hc)$ such that $g=ub$.  
\end{corollary}

\begin{proof}
By applying Theorem~\ref{A2} for $a=g^*g$ 
we get the operators $d\in\GL_{\Ig}$ and $r\in\Ig$ 
such that $g^*g=(\1+r^*)d(\1+r)$. 
Now denote $c=\1+r\in\GL_{\Ig}\cap\Alg\Pg$. 
Then $g^*g=c^*dc$, $d\ge0$, and all of the operators $g$, $c$, and $d$ 
are invertible, 
hence the operator $u:=g(d^{1/2}c)^{-1}$ is unitary. 
On the other hand, since $0\le d\in\GL_{\Ig}$, 
it is straightforward to prove that $d^{1/2}\in\GL_{\Ig}$, 
whence $u\in\U_{\Ig}(\Hc)$. 

In addition we have $b:=d^{1/2}c\in\GL_{\Ig}\cap\Alg\Pg$ 
and $g=ub$, and this completes the proof. 
\end{proof}

\begin{example}\label{A4}
\normalfont 
Theorem~\ref{A2} and Corollary~\ref{A3} apply in particular 
for the Schatten ideal $\Ig=\Sg_p(\Hc)$ if~$1<p<\infty$. 
\qed
\end{example}

\addtocontents{toc}{\SkipTocEntry}
\section*{Acknowledgment} 

The author wishes to thank 
Professor Hendrik Grundling and Professor Gary Weiss 
for kindly drawing his atention to some relevant references, 
as well as Professor Mihai \c Sabac for some useful remarks.

This work was partially supported by Grant
2-CEx06-11-34/25.07.06 of the Romanian Government.


\begin{thebibliography}{100000000}




\bibitem[Ara78]{Ara78}
J.~Arazy, 
Some remarks on interpolation theorems and the boundness of 
the triangular projection in unitary matrix spaces, 
\textit{Integral Equations Operator Theory} 
\textbf{1} (1978), no.~4, 453--495.

\bibitem[Arv67]{Arv67}
W.B.~Arveson, 
Analyticity in operator algebras,  
\textit{Amer. J. Math.} 
\textbf{89} (1967), 578--642. 

\bibitem[Arv75]{Arv75}
W.B.~Arveson, 
Interpolation problems in nest algebras,  
\textit{J. Functional Analysis} 
\textbf{20} (1975), no.~3, 208--233. 

\bibitem[Ba69]{Bal69}
V.K.~Balachandran, 
Simple $L^*$-algebras of classical type, 
\textit{Math. Ann.} 
\textbf{180} (1969), 205--219.

\bibitem[BD01]{BD01} 
V.~B\u alan, J.~Dorfmeister, 
Birkhoff decompositions and Iwasawa decompositions for loop groups, 
\textit{Tohoku Math. J.} 
\textbf{53} (2001), 593–-615. 

\bibitem[Be05]{Be01}
D.~Belti\c t\u a, 
On Banach-Lie algebras, spectral decompositions and complex polarizations. 
In: D. Ga\c spar, I.~Gohberg, D. Timotin, F.-H. Vasilescu, L. Zsido (eds.),
{\it Recent Advances in Operator Theory, Operator Algebras, 
and Their Applications.
XIXth International Conference on Operator Theory, Timisoara (Romania), 2002}.
Operator Theory: Advances and Applications, 153. Birkh\"auser Verlag
Basel,  
2005, pp.~13--38. 

\bibitem[Be06]{Bel06}
D.~Belti\c t\u a, 
\textit{Smooth Homogeneous Structures in Operator Theory}, 
Chapman \& Hall/CRC Monographs and Surveys in Pure
   and Applied Mathematics, 137. Chapman \& Hall/CRC Press, 
Boca Raton-London-New York-Singapore, 2006.

\bibitem[BP07]{BP05}
D.~Belti\c t\u a, B.~Prunaru,
Amenability, completely bounded projections,
dynamical systems and smooth orbits,
\textit{Integral Equations Operator Theory} 
(to appear). 
(See \textit{preprint} math.OA/0504313.)

\bibitem[BR05]{BR05}
D.~Belti\c t\u a, T.S.~Ratiu,
 Symplectic leaves in real Banach Lie-Poisson spaces, 
\textit{Geom. Funct. Anal.} 
\textbf{15} (2005), no.~4, 753--779. 

\bibitem[BR07]{BR06}
D.~Belti\c t\u a, T.S.~Ratiu, 
Geometric representation theory for unitary groups of operator algebras, 
\textit{Adv. Math.} \textbf{208} (2007), no.~1, 299--317.

\bibitem[BS01]{BS01}
D.~Belti\c t\u a, M.~\c Sabac,
{\it Lie Algebras of Bounded Operators},
Operator Theory: Advances and Applications, 120. Birkh\"auser Verlag,
Basel, 2001.

\bibitem[BFR93]{BFR93}
A.M.~Bloch, H.~Flaschka, T.S.~Ratiu, 
A Schur-Horn-Kostant convexity theorem for 
the diffeomorphism group of the annulus,
\textit{Invent. Math.} \textbf{113} (1993), no.~3, 511--529. 

\bibitem[Bo80]{Boy80}
R.P.~Boyer,
Representation theory of the Hilbert-Lie group
${\rm U}({\mathfrak H})_2$,
\textit{Duke Math. J.} \textbf{47} (1980), no.~2,
325--344.

\bibitem[Bo93]{Boy93}
R.P.~Boyer,
Representation theory of infinite-dimensional unitary
groups,
in:
\textit{Representation theory of groups and algebras},
 Contemp. Math., 145, Amer. Math. Soc., Providence, RI,
1993, pp.~381-391.

\bibitem[Ca85]{Ca85}
A.L.~Carey, 
Some homogeneous spaces and representations of 
the Hilbert Lie group ${\mathcal U}(H)_2$, 
\textit{Rev. Roumaine Math. Pures Appl.} 
\textbf{30} (1985), no.~7, 505--520. 

\bibitem[CG99]{CG99}
G.~Corach, J.E.~Gal\'e, 
On amenability and geometry of spaces of bounded representations, 
\textit{J. London Math. Soc. (2)} \textbf{59} (1999), no.~1, 311--329. 

\bibitem[Da88]{Dav88}
K.R.~Davidson, 
{\it Nest Algebras},
Pitman Research Notes in Mathematics Series, 191. 
Longman Scientific \& Technical, Harlow; 
copublished in the United States with 
John Wiley \& Sons, Inc., New York, 1988.

\bibitem[DPW02]{DPW02}
I.~Dimitrov, I.~Penkov, J.A.~Wolf, 
A Bott-Borel-Weil theory for direct limits of
algebraic groups, 
\textit{Amer. J. Math.} \textbf{124} (2002), 
no.~5, 955--998. 

\bibitem[DFWW04]{DFWW04}
K. Dykema, T. Figiel, G. Weiss, M. Wodzicki,
Commutator structure of operator ideals,
{\it Adv. Math.} 
{\bf 185}(2004),  no. 1, 1--79.

\bibitem[Er72]{Erd72}
J.A.~Erdos,
The triangular factorization of operators on Hilbert space,  
\textit{Indiana Univ. Math. J.} 
\textbf{22} (1972/73), 939--950.

\bibitem[Er78]{Erd78}
J.A.~Erdos, 
Triangular integration on symmetrically normed ideals, 
\textit{Indiana Univ. Math. J.} 
\textbf{27} (1978), no.~3, 401--408.

\bibitem[EL72]{EL72}
J.A.~Erdos, W.E.~Longstaff, 
The convergence of triangular integrals of operators on Hilbert space, 
\textit{Indiana Univ. Math. J.} 
\textbf{22} (1972/73), 929--938. 

\bibitem[Ga06]{Ga06}
J.E.~Gal\'e, 
Geometr\'\i a de \'orbitas de representaciones de grupos y \'algebras  
promediables, 
\textit{Rev. R. Acad. Cienc. Exactas F\'\i s. Qu\'\i m. Nat. Zaragoza (2)} 
(to appear).  

\bibitem[GK69]{GK69}
I.C.~Gohberg, M.G.~Kre\u\i n,
\textit{Introduction to the Theory of Linear Nonselfadjoint 
Operators}, 
Translations of Mathematical Monographs, vol.~18, 
American Mathematical Society, Providence, R.I., 1969.

\bibitem[GK70]{GK70}
I.C.~Gohberg, M.G.~Kre\u\i n, 
\textit{Theory and Applications of Volterra Operators in Hilbert Space}. 
Translated from the Russian by A.~Feinstein. 
Translations of Mathematical Monographs, Vol. 24, 
American Mathematical Society, Providence, R.I., 1970.

\bibitem[Gru05]{Gru05}
H.~Grundling, 
Generalising group algebras,  
\textit{J. London Math. Soc. (2)} 
\textbf{72} (2005), no.~3, 742--762.

\bibitem[Ha72]{dlH72}
P.~de la Harpe, 
\textit{Classical Banach-Lie Algebras and Banach-Lie Groups 
of Operators in Hilbert Space}, 
Lecture Notes in Mathematics 285, Springer-Verlag, 
Berlin-Heidelberg-New York, 1972. 

\bibitem[HK77]{HK77}
L.A. Harris, W. Kaup, 
Linear algebraic groups in infinite dimensions,
{\it Illinois J. Math.}
{\bf 21}(1977), 666--674.


\bibitem[He01]{Hel01}
S.~Helgason, 
\textit{Differential Geometry, Lie Groups, and Symmetric Spaces}. 
(Corrected reprint of the 1978 original). 
Graduate Studies in Mathematics, 34. 
American Mathematical Society, Providence, RI, 2001. 

\bibitem[Iw49]{Iwa49}
K.~Iwasawa, 
On some types of topological groups,  
\textit{Ann. of Math. (2)} 
\textbf{50} (1949), 507--558. 

\bibitem[KW02]{KW02}
V.~Kaftal, G.~Weiss, 
Traces, ideals, and arithmetic means,  
\textit{Proc. Natl. Acad. Sci. USA} 
\textbf{99} (2002), no.~11, 7356--7360.

\bibitem[KW06]{KW06}
V.~Kaftal, G.~Weiss, 
$B(H)$ lattices, density and arithmetic mean ideals, 
\textit{preprint}, 2006.

\bibitem[Ke04]{Kel99} 
P.~Kellersch, 
\textit{Eine Verallgemeinerung der Iwasawa Zerlegung in Loop Gruppen}, 
(Dissertation, TU-M\"{u}nchen, 1999), 
Differential Geometry - Dynamical Systems. Monogr.~4. 
Geometry Balkan Press, Bucharest, 2004.

\bibitem[Ki73]{Ki73}
A.A.~Kirillov, 
Representations of the infinite-dimensional unitary group. (Russian) 
\textit{Dokl. Akad. Nauk. SSSR} 
\textbf{212} (1973), 288--290. 

\bibitem[Kn96]{Kna96}
A.W.~Knapp, 
\textit{Lie Groups Beyond an Introduction}, 
Progress in Mathematics 140, Birkh\"auser-Verlag, 
Boston-Basel-Berlin, 1996.

\bibitem[Ko73]{Ko73}
B.~Kostant, 
On convexity, the Weyl group and the Iwasawa decomposition,  
\textit{Ann. Sci. \'Ecole Norm. Sup. (4)} \textbf{6} (1973), 413--455. 

\bibitem[Lan01]{Lan01}
S.~Lang,
{\it Fundamentals of Differential Geometry} (corrected second printing),
Graduate Texts in Mathematics, 191. Springer-Verlag,
New-York,
2001.

\bibitem[Lar85]{Lar85}
D.R.~Larson, 
Nest algebras and similarity transformations, 
\textit{Ann. of Math. (2)} \textbf{121} (1985), 
no.~3, 409--427.

\bibitem[LR91]{LR91}
J.-H.~Lu, T.S.~Ratiu, 
On the nonlinear convexity theorem of Kostant, 
\textit{J. Amer. Math. Soc.} \textbf{4} (1991), no.~2,
349--363. 

\bibitem[MSS88]{MSS88}
P.S.~Muhly, K.-S.~Saito, B.~Solel, 
Coordinates for triangular operator algebras, 
\textit{Ann. of Math. (2)}
\textbf{127} (1988), no.~2, 245--278.

\bibitem[NRW01]{NRW01}
L.~Natarajan, E.~Rodr\'\i guez-Carrington, J.A.~Wolf, 
The Bott-Borel-Weil theorem for direct limit groups,  
\textit{Trans. Amer. Math. Soc.} \textbf{353} (2001), 
no.~11, 4583--4622. 

\bibitem[Nee98]{Nee98}
K.-H.~Neeb, 
Holomorphic highest weight representations of
   infinite-dimensional complex classical groups, 
\textit{J. reine angew. Math.}
\textbf{497} (1998), 171--222.

\bibitem[Nee02a]{Ne02b}
K.-H.~Neeb, 
Classical Hilbert-Lie groups, their extensions 
and their homotopy groups, 
in: A. Strasburger, J. Hilgert, K.-H. Neeb, W. Wojyy\'nski (eds.),  
{\it Geometry and Analysis on Finite and Infinite-dimensional Lie Groups 
(B\c ed\l ewo, 2000)},
Banach Center Publ., vol. 55, Polish Acad. Sci. 
Warsaw, 
2002, pp.~87--151.

\bibitem[Nee02b]{Nee02c}
K.-H.~Neeb,
A Cartan-Hadamard theorem for Banach-Finsler manifolds,  
\textit{Geom. Dedicata} \textbf{95} (2002), 115--156. 


\bibitem[Nee04]{Nee04}
K.-H.~Neeb, 
Infinite-dimensional groups and their representations,
in: \textit{Lie Theory},
Progr. Math. 228, Birk\-h\"auser, Boston, MA, 2004,
pp.~213--328.

\bibitem[N\O 98]{NO98}
K.-H.~Neeb, B.~\O rsted, 
Unitary highest weight representations in Hilbert spaces 
of holomorphic functions on infinite-dimensional domains,  
\textit{J. Funct. Anal.} \textbf{156} (1998),
no.~1, 263--300.

\bibitem[Neu99]{Neu99}
A.~Neumann, 
An infinite-dimensional version of the Schur-Horn convexity theorem,  
\textit{J. Funct. Anal.} \textbf{161} (1999),
no.~2, 418--451.

\bibitem[Neu02]{Neu02}
A.~Neumann, 
An infinite dimensional version of the Kostant convexity theorem, 
\textit{J. Funct. Anal.} \textbf{189} (2002), no.~1,
80--131.

\bibitem[Ol78]{Ol78}
G.I.~Ol'\v sanski\u\i, 
Unitary representations of the infinite-dimensional classical 
groups $\text{U}(p,\infty )$,
$\text{SO}\sb{0}(p,\infty )$, $\text{Sp}(p,\infty )$, 
and of the corresponding motion groups. (Russian) 
\textit{Funktsional. Anal. i Prilozhen.} 
\textbf{12} (1978), no.~3, 32--44, 96. 

\bibitem[Ol88]{Ol88}
G.I.~Ol'\v sanski\u\i,
The method of holomorphic extensions in 
the theory of unitary representations of
infinite-dimensional classical groups. (Russian) 
\textit{Funktsional. Anal. i Prilozhen.} 
\textbf{22} (1988), no.~4, 23--37, 96.  

\bibitem[Pa88]{Pat88}
A.L.T.~Paterson, 
\textit{Amenability}, 
Mathematical Surveys and Monographs, 29. American Mathematical Society,
Providence, RI, 1988.

\bibitem[Pic90]{Pic90}
D.~Pickrell, 
Separable representations for automorphism groups of 
infinite symmetric spaces, 
\textit{J. Funct. Anal.} 
\textbf{90} (1990), no.~1, 1--26.

\bibitem[Pit88]{Pit88}
D.R.~Pitts, 
Factorization problems for nests: factorization methods and 
characterizations of the universal factorization property, 
\textit{J. Funct. Anal.} \textbf{79} (1988), no.~1, 57--90. 

\bibitem[Po86]{Pow86}
S.C.~Power, 
Factorization in analytic operator algebras,  
\textit{J. Funct. Anal.} 
\textbf{67} (1986), no.~3, 413--432.


\bibitem[Se57]{Se57}
I.E.~Segal, 
The structure of a class of representations of 
the unitary group on a Hilbert space, 
\textit{Proc. Amer. Math. Soc.} 
\textbf{8} (1957), 197--203.

\bibitem[SV75]{SV75}
\c S.~Str\u atil\u a, D.~Voiculescu, 
\textit{Representations of AF-algebras and of 
the group $\text{U}(\infty )$}, 
Lecture Notes in
Mathematics, Vol.~486. Springer-Verlag, Berlin-New York, 1975.

\bibitem[Tu05]{Tum05}
A.B.~Tumpach, 
\textit{Vari\'et\'es K\"ahl\'eriennes et
Hyperk\"ahl\'eriennes de Dimension Infinie}, 
Ph.D Thesis, 
\'Ecole Polytechnique, Paris, 2005.

\bibitem[Tu06]{Tum06} 
A.B.~Tumpach, 
Mostow decomposition theorems
for $L^*$-groups and applications to affine coadjoint orbits and
stable manifolds, 
\textit{preprint} math-ph/0605039.

\bibitem[Up85]{Up85} 
H.~Upmeier, 
\textit{Symmetric Banach Manifolds and Jordan $C^*$-Algebras}, 
North-Holland Math. Stud. 104, Notas de Matem\'atica 96, 
North-Holland, Amsterdam, 1985. 


\bibitem[We05]{We05}
G.~Weiss, $B(H)$-commutators: a historical survey. 
In: D. Ga\c spar, I. Gohberg, D. Timotin, F.-H. Vasilescu, L. Zsido (eds.),
{\it Recent Advances in Operator Theory, Operator Algebras, 
and Their Applications.
XIXth International Conference on Operator Theory, Timisoara (Romania), 2002}.
Operator Theory: Advances and Applications, 153. Birkh\"auser Verlag
Basel,  
2005, pp.~307--320.


\bibitem[Wo05]{Wo04}
J.A.~Wolf,
Principal series representations of direct limit
groups,
\textit{Compos. Math.} \textbf{141} (2005), no.~6, 1504--1530.


\end{thebibliography}
\end{document}